\documentclass[12pt,draft]{article}

\usepackage{color}
\usepackage{amsmath}
\usepackage{amssymb}
\usepackage{amsbsy}

\usepackage{amsfonts}
\usepackage{tikz}

\usepackage[margin=1in]{geometry}

\newtheorem{lem}{Lemma}[section]%
\newtheorem{thm}[lem]{Theorem}%
\newtheorem{ex}[lem]{Example}%
\newtheorem{defi}[lem]{Definition}%
\newtheorem{prob}{Problem}%
\newtheorem{prop}[lem]{Proposition}%
\newtheorem{rem}[lem]{Remark}%

\def\a{\alpha}    
 \def\s{\sigma}

\def\G{\Gamma}

 \def\lg{\langle} \def\rg{\rangle}
\def\nd{\mathrel{\bigm|\kern-.7em/}}

\def\f{\noindent}

\def\PSL{\hbox{\rm PSL}}
\def\SL{\hbox{\rm SL}} \def\PGL{\hbox{\rm PGL}}\def\GL{\hbox{\rm GL}}

\DeclareMathOperator{\Aut}{Aut}

\def\soc{\hbox{\rm soc}}

\def\Cay{\hbox{\rm Cay}}

\DeclareMathOperator{\Sym}{Sym}

\def\mz{{\mathbb Z}}
\def\K{{\rm\bf K}}
\def\C{{\bf C}}

\def\I{\mathcal{G}}
\def\II{\widetilde{\mathcal{G}}}

\def\demo{\noindent{\bf Proof}\hskip10pt}

\def\qed{\hskip10pt $\Box$\vspace{3mm}}

\newcommand{\F}{\mathbb{F}}

\tikzset{every picture/.style={line width=0.75pt}} 

\begin{document}
\title{A characterisation of edge-affine $2$-arc-transitive covers of $\K_{2^n,2^n}$}
\author
{
Daniel R. Hawtin\\
{\small Faculty of Mathematics, The University of Rijeka}\\[-4pt]
{\small Rijeka, 51000, Croatia}\\[-2pt]
{\small{Email}:\ dan.hawtin@gmail.com}\\[+6pt]
Cheryl E. Praeger\\
{\small Department of Mathematics and Statistics, The University of Western Australia}\\[-4pt]
{\small Crawley, WA 6907, Australia}\\[-2pt]
{\small{Email}:\ cheryl.praeger@uwa.edu.au}\\[+6pt]
Jin-Xin Zhou  \\
{\small School of Mathematics and Statistics, Beijing Jiaotong University}\\[-4pt]
{\small Beijing 100044, P.R. China}\\[-2pt]
{\small{Email}:\ jxzhou@bjtu.edu.cn}
}

\date{}
\maketitle
\begin{abstract}
We introduce the notion of an \emph{$n$-dimensional mixed dihedral group}, a general class of groups for which we give a graph theoretic characterisation.
In particular, if $H$ is an $n$-dimensional mixed dihedral group then the we construct an edge-transitive Cayley graph $\Gamma$ of $H$ such that the clique graph $\Sigma$ of $\Gamma$ is a $2$-arc-transitive normal cover of $\K_{2^n,2^n}$, with a subgroup of $\Aut(\Sigma)$ inducing a particular \emph{edge-affine} action on $\K_{2^n,2^n}$. Conversely, we prove that if $\Sigma$ is a $2$-arc-transitive normal cover of $\K_{2^n,2^n}$, with a subgroup of $\Aut(\Sigma)$ inducing an \emph{edge-affine} action on $\K_{2^n,2^n}$, then the line graph $\Gamma$ of $\Sigma$ is a Cayley graph of an $n$-dimensional mixed dihedral group.




Furthermore, we give an explicit construction of a family of $n$-dimensional mixed dihedral groups. This family addresses a problem proposed by Li concerning normal covers of prime power order of the `basic' $2$-arc-transitive graphs. In particular, we construct, for each $n\geq 2$, a $2$-arc-transitive normal cover of $2$-power order of the `basic' graph $\K_{2^n,2^n}$.

\bigskip
\noindent{\bf Key words:} $2$-arc-transitive, normal cover, Cayley graph, edge-transitive\\
\noindent{\bf 2000 Mathematics subject classification:} 05C38, 20B25
\end{abstract}

\section{Introduction}\label{s:intro}




All graphs we consider are finite, connected, simple and undirected. We denote the vertex set of a graph $\Gamma$ by $V(\Gamma)$ and its edge set by $E(\Gamma)$, and call $|V(\G)|$ the \emph{order} of $\G$. For a subgroup $G$ of the automorphism group $\Aut(\G)$, a graph $\G$ is said to be  \emph{$(G,2)$-arc-transitive} if $G$ acts transitively on the set of all \emph{$2$-arcs} of $\G$ (that is, the set of all triples $(u,v,w)$ such that $u,v,w\in V(\G)$, $u\neq w$ and $\{u,v\},\{v,w\}\in E(\G)$); and $\Gamma$ is \emph{$2$-arc-transitive} if such a group exists.
A graph $\G$ is called a \emph{normal cover}  of a graph $\Sigma$ if $\G$ and $\Sigma$ have the same valency, and there exists a subgroup chain $N\lhd G\leq\Aut(\G)$ such that $G$ is transitive on $V(\G)$ and $\Sigma$ is the quotient of $\G$ by the set of $N$-orbits in $V(\G)$. We also say that $\G$ is an \emph{$N$-normal cover} of $\Sigma$
to specify the role of $N$, and we call $\Sigma$ the  \emph{$N$-normal quotient} of $\G$, often denoted $\G_N$.

In \cite{L-BLMS-2001}, Li proved that every $2$-arc-transitive graph of prime power order is a normal cover of one of the following graphs (which we regard as `basic graphs'): $\K_{2^n,2^n}$ (the complete bipartite graph), $\K_{p^m}$ (the complete graph), $\K_{2^n,2^n}-2^n\K_2$ (the graph obtained by deleting a $1$-factor from $\K_{2^n,2^n}$) or a primitive or biprimitive affine graph (by which, see \cite[p.130]{L-BLMS-2001}, Li meant the graphs  in the classification by Ivanov and the second author in \cite[Table 1]{IP}). Li was ``inclined to think that non-basic 2-arc-transitive graphs of prime power order would be rare and hard to construct", see \cite[pp.130--131]{L-BLMS-2001}, and posed the following problem.

\begin{prob}{\rm\cite[Problem]{L-BLMS-2001}}\label{prob-2}
Construct and characterize the normal covers of prime power order, of the basic $2$-arc-transitive graphs of prime power order.
\end{prob}

As mentioned in \cite{L-BLMS-2001}, several non-trivial normal covers of order a $2$-power of primitive affine graphs appear in \cite[Theorem 1.4]{Ivanov-Praeger}. These are certain girth $5$ normal covers of the Armanios--Wells graph of order $32$, of which only the smallest one (of order $2^{20}$) has been confirmned to exist --  by a coset enumeration performed by Leonard Soicher, \cite[Section 8]{IP}.

In this paper, we give a partial answer to Problem~\ref{prob-2} by characterising a certain family of $2$-arc-transitive normal covers of $\K_{2^n,2^n}$, and by explicitly constructing an infinite family of such graphs, which are all non-basic and have $2$-power order.

In addition, the family of graphs we construct provides an affirmative answer to Problem~\ref{prob-3} below, which was posed by Chen et al.~in \cite{Chen-Jin-Li2019}. A non-complete graph $\G$ of diameter $d\geq 2$ is called {\em $2$-distance-transitive} if $\Aut(\G)$ acts transitively, for each $i=0,1,2$, on the set of all pairs $(u,v)$, where $u,v\in V(\G)$ and $d(u,v)=i$, and {\em distance-transitive} if $\Aut(\G)$ acts transitively, for each $i=0,1,\ldots,d$, on the set of all pairs $(u,v)$, where $u,v\in V(\G)$ and $d(u,v)=i$. (See \cite{Chen-Jin-Li2019,CJin-S-JGT2017,JHL} for more information about $2$-distance-transitive graphs.) Also, a Cayley graph $\Cay(G,S)$ is said to be {\em normal} if $G$ is normal in $\Aut(\Cay(G,S))$ (see Section~\ref{sub:cay}). Note that Problem~\ref{prob-3} has been answered in the affirmative in \cite{huang2022distance}, where an infinite family of such graphs was constructed. The graphs we construct here are different to those in \cite{huang2022distance}.



\begin{prob}{\rm \cite[Question~1.2]{Chen-Jin-Li2019}}\label{prob-3}
Is there a normal Cayley graph that is $2$-distance-transitive, but is neither distance-transitive nor $2$-arc-transitive?
\end{prob}

\subsection{Statement of main results}

In order to state our main results we introduce in Definition~\ref{mixed-dihedral-group} a class of groups, and for each such group we introduce  in Definition~\ref{mixed-dihedrant} two associated graphs. Note that the \emph{derived subgroup} $G'$  of a group $G$  is the subgroup generated by all the {\em commutators}, that is, elements $[a,b]=a^{-1}b^{-1}ab$, for $a, b\in G$.

\begin{defi}\label{mixed-dihedral-group}
{\em Let $n\geq1$. If $H$ is a finite group with subgroups $X,Y$ such that $X\cong Y\cong C_2^n$, $H=\lg X\cup Y\rg$ and $H/H'\cong C_2^{2n}$, then we call $H$ an {\em $n$-dimensional mixed dihedral group relative to $X$ and $Y$}.}
\end{defi}

For $H, X, Y$ as in Definition~\ref{mixed-dihedral-group},  we show in Lemma~\ref{prop-of-mixed-dihedral} that the map $\phi: H\to H/H'$ given by $h\mapsto hH'$ is a homomorphism such that $H/H'=\phi(X)\times\phi(Y)$.
We note that, for each $n\geq 1$, there exist infinitely many $n$-dimensional mixed dihedral groups. For example, let $m_1,m_2,\cdots, m_n$ be positive even integers, and let
\[\begin{array}{l}
A_i=\lg x_i,y_i\mid x_i^2=y_i^2=(x_iy_i)^{m_i}=1\rg\cong D_{2m_i}\ {\rm with}\ i=1,2,\ldots, n.\\
\end{array}\]
It is easy to see that $A_1\times A_2\times\cdots\times A_n\cong\prod_{i=1}^n D_{2m_i}$ is an $n$-dimensional mixed dihedral group relative to $\lg x_1,\ldots,x_n\rg$ and $\lg y_1,\ldots,y_n\rg$. These natural examples suggest the name for this family of groups but, as we shall see, there are many $n$-dimensional mixed dihedral groups apart from direct products of dihedral groups.

We now introduce a pair of graphs for each mixed dihedral group. See Subsection~\ref{sub:cay} for definitions related to the Cayley graph $\Cay(G,S)$.

\begin{defi}\label{mixed-dihedrant}
{\em Let $n\geq 2$, and let $H$ be an $n$-dimensional mixed dihedral group relative to $X$ and $Y$, as in Definition~\ref{mixed-dihedral-group}. Define the graphs $C(H,X,Y)$ and $\Sigma(H,X,Y)$ as follows:
\begin{equation}\label{eq-1}
C(H,X,Y)=\Cay(H,S(X,Y)),\ {\rm with}\ S(X,Y)=(X\cup Y)\setminus\{1\};
\end{equation}
%
and for $\Sigma=\Sigma(H,X,Y)$,
\begin{equation}\label{eq-2}
\begin{array}{l}
V(\Sigma)=\{Xh : h\in H\}\cup \{Yh: h\in H\},\\
E(\Sigma)=\{\{Xh, Yg\}:  h,g\in H\ \mbox{and}\ Xh\cap Yg\neq\emptyset\}.
\end{array}
\end{equation}}
\end{defi}

\begin{rem}\label{remark-sigma-graph}
{\rm For $\Sigma=\Sigma(H, X, Y)$ in Definition~\ref{mixed-dihedrant},  the set $\Sigma(Xh)$ of vertices adjacent to the vertex $Xh$ is $\{ Yxh : x\in X\}$, and similarly, $\Sigma(Yg)=\{ Xyg : y\in Y\}$. We show in Lemma~\ref{prop-of-mixed-dihedral}~(2) that $X\cap Y=1$, and this implies that $|\Sigma(Xh)|=|\Sigma(Yg)|=2^n$ for all $h,g\in H$, and hence that $\Sigma$ is a regular bipartite graph of valency $2^n$.}
\end{rem}

In this paper we are interested in a particular class of normal covers of $\K_{2^n,2^n}$. The graph $\Sigma=\K_{2^n,2^n}$ is said to be {\em $H$-edge-affine} if $H\leq S_{2^n}\wr S_2=\Aut(\Sigma)$ and $H$ has a normal subgroup $C_{2}^n\times C_{2}^n$ that is intransitive on vertices and regular on edges. By Proposition~\ref{basic-auto}, a $(G,2)$-arc-transitive $N$-normal cover $\G$ of $\K_{2^n,2^n}$ comes in four `flavours', and one of these is the case where the normal quotient $\G_N$ is $G/N$-edge-affine. We prove the following, further motivating our main result.


\begin{thm}\label{is-affine}
Let $n\geq 2$ and let $N\lhd\, G\leq\Aut(\G)$, for a graph $\G$, such that $\G$ is $(G,2)$-arc-transitive and is an $N$-normal cover of $\K_{2^n,2^n}$. Then either $\G$ is a Cayley graph, or $n\geq4$ and $\G_N$ is $G/N$-edge-affine.
\end{thm}


We now discuss our main result Theorem~\ref{iff-characterization} which characterises $(G,2)$-arc-transitive $N$-normal covers $\G$ of $\K_{2^n,2^n}$ such that $\G_N$ is $G/N$-edge-affine. In particular Theorem~\ref{iff-characterization}  answers Problem~\ref{prob-2} for the basic graphs $\K_{2^n,2^n}$. It does more than this since we do not require the group $N$ to be a $2$-group, see Example~\ref{ex} for examples.
The {\em clique graph} $C(\G)$ of $\G$ is the graph with vertices the maximal cliques of $\G$ such that two different maximal cliques are adjacent if and only if they share at least one common vertex.
We denote the derived subgroup of a group $H$ by $H'$.


\begin{ex}\label{ex}
{\rm
 An infinite family of examples satisfying the conditions in Theorem~\ref{iff-characterization}, but for which $N$ is not a $2$-group, can be obtained from a construction of Poto\v cnik and Spiga in \cite{PS2019}: Let $G\leq\Aut(\K_{2^n,2^n})$ be $2$-arc-transitive on $\K_{2^n,2^n}$ with $\soc(G)\cong C_2^n\times C_2^n$ acting intransitively on the vertices and regularly on the edges of $\K_{2^n,2^n}$. By \cite[Theorem 6]{PS2019}, for each odd prime $p$ there exists a $q$-fold regular cover $\G$ of $\K_{2^n,2^n}$, for some power $q$ of $p$, such that the maximal lifted group of automorphisms of $\G$ induces precisely $G$ on $\K_{2^n,2^n}$. The lifted group has the form $P\rtimes G$, with $P$ a group of order $q$,  $P\rtimes G$ acts as a 2-arc-transitive group on  $\G$, and the lift of  $\soc(G)$ is $X:=P\rtimes\soc(G)$ acting regularly on the edges of  $\G$. For an edge $\{u,v\}$ of $\G$, we have $X_u\cong X_v\cong C_2^n$, $X=\lg X_u,  X_v\rg$ and $X/X'\cong C_2^{2n}.$
}
\end{ex}


\begin{thm}\label{iff-characterization}
Let $n\geq 2$, and let $\Sigma$ be a graph,  $G\leq \Aut(\Sigma)$, and $N\lhd G$. Then the following are equivalent.
 \begin{enumerate}
    \item[{\rm (a)}]  $\Sigma$ is a $(G,2)$-arc-transitive $N$-normal cover of $\K_{2^n,2^n}$ which is $G/N$-edge-affine;
    \item[{\rm (b)}] $G$ has a normal subgroup $H$ with $H'=N$ such that $H$ is an $n$-dimensional mixed dihedral group relative to $X$ and $Y$, the line graph of $\Sigma$ is $C(H,X,Y)$, and $C(H,X,Y)$ is $G$-edge-transitive.
 \end{enumerate}
\end{thm}

We now define  an infinite family of mixed dihedral groups. 

 \medskip

\begin{defi}\label{constructionI}
Let $n\geq 2$, and let $\II(n)=\lg x_1,\ldots,x_n,y_1,\ldots,y_n\rg$ be a finite $2$-group
with the following defining relations, where $1\leq i,j,k\leq n$:
\[\begin{array}{l}
x_i^2=y_i^2=1, [x_i,x_j]=[y_i,y_j]=1,[x_i,y_j]^2=1, [[x_i,y_j],x_k]=[[x_i,y_j],y_k]=1,
\end{array}\]
Let $X=\lg x_1,\ldots,x_n\rg$ and $Y=\lg y_1,\ldots,y_n\rg$.
\end{defi}


In Theorem~\ref{2-dist-tran} we show that the graph constructions in Definition~\ref{mixed-dihedrant} using the groups in Definition~\ref{constructionI} provide an infinite family of examples answering Problem~\ref{prob-3} in the affirmative, namely the graphs $C(\II(n),X,Y)$ have all the properties  sought in  Problem~\ref{prob-3}. Also, for each $n\geq2$, the graph $\Sigma(\II(n),X,Y)$ has all the properties
 required in  Problem~\ref{prob-2}   for the basic $2$-arc-transitive graph $\K_{2^n,2^n}$. A graph $\G$ is said to be \emph{$2$-geodesic-transitive} if $\G$ is vertex-transitive and $\Aut(\G)$ acts transitively on the set of all triples $(u,v,w)$ where $u,v,w\in V(\G)$ with $\{u,v\},\{v,w\}\in E(\G)$ and $\{u,w\}\notin E(\G)$. Note that, $2$-geodesic-transitivity implies $2$-distance-transitivity.


\begin{thm}\label{2-dist-tran}
Let $n\geq2$, let $\II(n), X, Y$ be as in Definition~\ref{constructionI},  let $\Gamma=C(\II(n),X,Y)$ and $\Sigma=\Sigma(\II(n),X,Y)$ be as in Definition~\emph{\ref{mixed-dihedrant}}, and let $N=\langle [x_i,y_j]\mid 1\leq i,j\leq n \rangle$. Then
\begin{enumerate}
    \item[{\rm (1)}] $\II(n)$ is an $n$-dimensional mixed dihedral group relative to $X$ and $Y$, with derived subgroup $\II(n)'=N \cong C_2^{n^2}$;
    \item[{\rm (2)}] $\Sigma$ is a $2$-arc-transitive $N$-normal cover of $\K_{2^n,2^n}$, of order $2^{n^2+n-1}$, which is $\II(n)/N$-edge-affine;   
    \item[{\rm (3)}] $\Gamma$ is a $2$-geodesic-transitive normal Cayley graph; moreover, $\Gamma$ is $2$-distance-transitive, but it is neither distance-transitive nor $2$-arc-transitive.
\end{enumerate}
\end{thm}

\begin{rem}\label{rem:II}
{\rm In Definition~\ref{constructionI} we give an explicit definition of a group $\I(n)$ in terms of mappings of vector spaces, and in Lemma~\ref{CCpresentation} we prove that this group $\I(n)$ is isomorphic to the group $\II(n)$ in Definition~\ref{constructionI}.
This concrete representation of $\II(n)$ is useful for deriving various of its properties. Moreover, in Lemma~\ref{fullAutomorphism} we determine the full automorphism group $A$ of the graph $C(\II(n),X,Y)$, prove that it is equal to the full automorphism group of $\Sigma(\II(n),X,Y)$ and prove that $\K_{2^n,2^n}$ is $A/N$-edge-affine, a stronger statement than Theorem~\ref{2-dist-tran}~(2). Knowing the full automorphism group of $C(\II(n),X,Y)$ and $\Sigma(\II(n),X,Y)$ turns out to be essential for the proof of Theorem~\ref{2-dist-tran}. Figure~\ref{distdiaggamm22} gives the distance diagram of $C(\II(2),X,Y)$. }
\end{rem}

\begin{figure}
 \[
\begin{tikzpicture}[x=0.75pt,y=0.75pt,yscale=-1,xscale=1]

\draw   (85,135) .. controls (85,126.72) and (91.72,120) .. (100,120) .. controls (108.28,120) and (115,126.72) .. (115,135) .. controls (115,143.28) and (108.28,150) .. (100,150) .. controls (91.72,150) and (85,143.28) .. (85,135) -- cycle ;
\draw   (435,195) .. controls (435,186.72) and (441.72,180) .. (450,180) .. controls (458.28,180) and (465,186.72) .. (465,195) .. controls (465,203.28) and (458.28,210) .. (450,210) .. controls (441.72,210) and (435,203.28) .. (435,195) -- cycle ;
\draw   (435,135) .. controls (435,126.72) and (441.72,120) .. (450,120) .. controls (458.28,120) and (465,126.72) .. (465,135) .. controls (465,143.28) and (458.28,150) .. (450,150) .. controls (441.72,150) and (435,143.28) .. (435,135) -- cycle ;
\draw   (365,195) .. controls (365,186.72) and (371.72,180) .. (380,180) .. controls (388.28,180) and (395,186.72) .. (395,195) .. controls (395,203.28) and (388.28,210) .. (380,210) .. controls (371.72,210) and (365,203.28) .. (365,195) -- cycle ;
\draw   (365,75) .. controls (365,66.72) and (371.72,60) .. (380,60) .. controls (388.28,60) and (395,66.72) .. (395,75) .. controls (395,83.28) and (388.28,90) .. (380,90) .. controls (371.72,90) and (365,83.28) .. (365,75) -- cycle ;
\draw   (365,135) .. controls (365,126.72) and (371.72,120) .. (380,120) .. controls (388.28,120) and (395,126.72) .. (395,135) .. controls (395,143.28) and (388.28,150) .. (380,150) .. controls (371.72,150) and (365,143.28) .. (365,135) -- cycle ;
\draw   (295,165) .. controls (295,156.72) and (301.72,150) .. (310,150) .. controls (318.28,150) and (325,156.72) .. (325,165) .. controls (325,173.28) and (318.28,180) .. (310,180) .. controls (301.72,180) and (295,173.28) .. (295,165) -- cycle ;
\draw   (295,105) .. controls (295,96.72) and (301.72,90) .. (310,90) .. controls (318.28,90) and (325,96.72) .. (325,105) .. controls (325,113.28) and (318.28,120) .. (310,120) .. controls (301.72,120) and (295,113.28) .. (295,105) -- cycle ;
\draw   (225,135) .. controls (225,126.72) and (231.72,120) .. (240,120) .. controls (248.28,120) and (255,126.72) .. (255,135) .. controls (255,143.28) and (248.28,150) .. (240,150) .. controls (231.72,150) and (225,143.28) .. (225,135) -- cycle ;
\draw   (155,135) .. controls (155,126.72) and (161.72,120) .. (170,120) .. controls (178.28,120) and (185,126.72) .. (185,135) .. controls (185,143.28) and (178.28,150) .. (170,150) .. controls (161.72,150) and (155,143.28) .. (155,135) -- cycle ;
\draw   (505,135) .. controls (505,126.72) and (511.72,120) .. (520,120) .. controls (528.28,120) and (535,126.72) .. (535,135) .. controls (535,143.28) and (528.28,150) .. (520,150) .. controls (511.72,150) and (505,143.28) .. (505,135) -- cycle ;
\draw    (115,135) -- (155,135) ;
\draw    (185,135) -- (225,135) ;
\draw    (255,130) -- (295,110) ;
\draw    (255,140) -- (295,160) ;
\draw    (325,170) -- (365,190) ;
\draw    (325,110) -- (365,130) ;
\draw    (325,100) -- (365,80) ;
\draw    (325,160) -- (365,140) ;
\draw    (395,195) -- (435,195) ;
\draw    (460,185) -- (510,145) ;
\draw    (310,120) -- (310,150) ;
\draw    (380,90) -- (380,120) ;
\draw    (380,150) -- (380,180) ;
\draw    (390,185) -- (440,145) ;
\draw    (450,150) -- (450,180) ;

\draw (96,127.4) node [anchor=north west][inner sep=0.75pt]  [font=\small]  {$1$};
\draw (441,187.4) node [anchor=north west][inner sep=0.75pt]  [font=\small]  {$36$};
\draw (441,127.4) node [anchor=north west][inner sep=0.75pt]  [font=\small]  {$18$};
\draw (371,187.4) node [anchor=north west][inner sep=0.75pt]  [font=\small]  {$72$};
\draw (371,127.4) node [anchor=north west][inner sep=0.75pt]  [font=\small]  {$36$};
\draw (376,67.4) node [anchor=north west][inner sep=0.75pt]  [font=\small]  {$9$};
\draw (301,157.4) node [anchor=north west][inner sep=0.75pt]  [font=\small]  {$36$};
\draw (301,97.4) node [anchor=north west][inner sep=0.75pt]  [font=\small]  {$18$};
\draw (231,127.4) node [anchor=north west][inner sep=0.75pt]  [font=\small]  {$18$};
\draw (166,127.4) node [anchor=north west][inner sep=0.75pt]  [font=\small]  {$6$};
\draw (516,127.4) node [anchor=north west][inner sep=0.75pt]  [font=\small]  {$6$};
\draw (117,116.4) node [anchor=north west][inner sep=0.75pt]  [font=\footnotesize]  {$6$};
\draw (146,116.4) node [anchor=north west][inner sep=0.75pt]  [font=\footnotesize]  {$1$};
\draw (166,102.4) node [anchor=north west][inner sep=0.75pt]  [font=\footnotesize]  {$2$};
\draw (186,116.4) node [anchor=north west][inner sep=0.75pt]  [font=\footnotesize]  {$3$};
\draw (216,116.4) node [anchor=north west][inner sep=0.75pt]  [font=\footnotesize]  {$1$};
\draw (236,102.4) node [anchor=north west][inner sep=0.75pt]  [font=\footnotesize]  {$2$};
\draw (256,111.4) node [anchor=north west][inner sep=0.75pt]  [font=\footnotesize]  {$1$};
\draw (256,146.4) node [anchor=north west][inner sep=0.75pt]  [font=\footnotesize]  {$2$};
\draw (502,151.4) node [anchor=north west][inner sep=0.75pt]  [font=\footnotesize]  {$6$};
\draw (467,182.4) node [anchor=north west][inner sep=0.75pt]  [font=\footnotesize]  {$1$};
\draw (426,196.4) node [anchor=north west][inner sep=0.75pt]  [font=\footnotesize]  {$2$};
\draw (426,137.4) node [anchor=north west][inner sep=0.75pt]  [font=\footnotesize]  {$4$};
\draw (391,166.4) node [anchor=north west][inner sep=0.75pt]  [font=\footnotesize]  {$1$};
\draw (396,131.4) node [anchor=north west][inner sep=0.75pt]  [font=\footnotesize]  {$1$};
\draw (367,151.4) node [anchor=north west][inner sep=0.75pt]  [font=\footnotesize]  {$2$};
\draw (396,196.4) node [anchor=north west][inner sep=0.75pt]  [font=\footnotesize]  {$1$};
\draw (352,187.4) node [anchor=north west][inner sep=0.75pt]  [font=\footnotesize]  {$1$};
\draw (351,131.4) node [anchor=north west][inner sep=0.75pt]  [font=\footnotesize]  {$1$};
\draw (357,111.4) node [anchor=north west][inner sep=0.75pt]  [font=\footnotesize]  {$1$};
\draw (351,67.4) node [anchor=north west][inner sep=0.75pt]  [font=\footnotesize]  {$2$};
\draw (446,212.4) node [anchor=north west][inner sep=0.75pt]  [font=\footnotesize]  {$2$};
\draw (376,212.4) node [anchor=north west][inner sep=0.75pt]  [font=\footnotesize]  {$2$};
\draw (371,166.4) node [anchor=north west][inner sep=0.75pt]  [font=\footnotesize]  {$1$};
\draw (382,91.4) node [anchor=north west][inner sep=0.75pt]  [font=\footnotesize]  {$4$};
\draw (326,176.4) node [anchor=north west][inner sep=0.75pt]  [font=\footnotesize]  {$2$};
\draw (326,142.4) node [anchor=north west][inner sep=0.75pt]  [font=\footnotesize]  {$1$};
\draw (326,116.4) node [anchor=north west][inner sep=0.75pt]  [font=\footnotesize]  {$2$};
\draw (326,81.4) node [anchor=north west][inner sep=0.75pt]  [font=\footnotesize]  {$1$};
\draw (302,182.4) node [anchor=north west][inner sep=0.75pt]  [font=\footnotesize]  {$1$};
\draw (282,161.4) node [anchor=north west][inner sep=0.75pt]  [font=\footnotesize]  {$1$};
\draw (282,96.4) node [anchor=north west][inner sep=0.75pt]  [font=\footnotesize]  {$1$};
\draw (297,121.4) node [anchor=north west][inner sep=0.75pt]  [font=\footnotesize]  {$2$};
\draw (297,137.4) node [anchor=north west][inner sep=0.75pt]  [font=\footnotesize]  {$1$};
\draw (381,107.4) node [anchor=north west][inner sep=0.75pt]  [font=\footnotesize]  {$1$};
\draw (441,167.4) node [anchor=north west][inner sep=0.75pt]  [font=\footnotesize]  {$1$};
\draw (452,151.4) node [anchor=north west][inner sep=0.75pt]  [font=\footnotesize]  {$2$};

\end{tikzpicture}
 \]
 \caption{Distance diagram of the graph $C(\II(2),X,Y)$, as in Theorem~\ref{2-dist-tran}, for the orbits of the stabiliser of a vertex inside the full automorphism group. Computations performed in GAP \cite{GAP4}.  
 }
 \label{distdiaggamm22}
\end{figure}

This paper is organised as follows. In Section~\ref{prelimSect}, we outline the notation used in the paper and give several preliminary results. In Section~3, we prove Theorem~\ref{is-affine}; in Section~4 we prove Theorem~\ref{iff-characterization}; and in Section~5 we use Theorem~\ref{iff-characterization} to prove Theorem~\ref{2-dist-tran}.

\section{Preliminaries}\label{prelimSect}

In this section, we introduce the notation and concepts we require for graphs and their symmetry properties, and we prove some preliminary results.

\subsection{Notation and concepts for graphs and groups}

We begin with various notions we will meet for graphs.

\subsection{Concepts for graphs}
Let $\G$ be a graph. As in Section~\ref{s:intro}, $V(\Gamma)$, $E(\Gamma)$ and $\Aut(\Gamma)$  denote its vertex set, edge set, and full automorphism group, respectively. Let $d(\G)$ be the diameter of $\Gamma$ (the maximum distance between vertices).  For $v\in V(\G)$ and $1\leq i\leq d(\G)$, let $\G_i(v)$ denote the set of vertices at distance $i$ from $v$; we often write $\G(v)=\G_1(v)$.  A graph is said to be {\em regular} if there exists an integer $k$ such that $|\G(v)|=k$ for all vertices $v\in V(\G)$.

A \emph{clique} of a graph $\G$ is a subset  $U\subseteq
V(\G)$ such that every pair of vertices in $U$ forms an
edge of $\G$. A clique $U$ is \emph{maximal} if no subset
of $V(\G)$  properly containing $U$ is a clique. The \emph{clique graph}  of $\G$ is defined as the graph $\Sigma(\G)$ with vertices the maximal cliques of $\G$ such that two distinct maximal cliques are adjacent in $\Sigma(\G)$ if and only if their intersection is non-empty. Similarly the \emph{line graph}  of $\G$ is defined as the graph $\mathcal{L}(\G)$ with vertex set $E(\G)$ such that two distinct edges $e,e'\in E(\G)$ are adjacent in $\mathcal{L}(\G)$ if and only if $e\cap e'\ne\emptyset$.

A graph $\G$ is bipartite if $E(\G)\ne\emptyset$ and $V(\G)$ is of the form $\Delta\cup\Delta'$ such that each edge consists of one vertex from $\Delta$ and one vertex from $\Delta'$. If $\G$ is connected then this vertex partition is uniquely determined and the two parts $\Delta, \Delta'$ are often called the \emph{biparts} of $\G$.

\subsubsection{Symmetry concepts for graphs}

Let $G\leq\Aut(\G)$. For $v\in V(\G)$, let $G_v=\{g\in G\ :\ v^g=v\}$, the stabiliser of $v$ in $G$. We say that $\G$ is {\em $G$-vertex-transitive\/} or {\em $G$-edge-transitive\/} if $G$ is transitive on $V(\G)$ or $E(\G)$, respectively.
When $G=\Aut(\G)$, a $G$-vertex-transitive or $G$-edge-transitive graph $\G$ is simply called {\em vertex-transitive\/} or {\em edge-transitive\/}, respectively. A regular graph $\G$ is said to be {\em $G$-locally primitive\/} or {\em locally $(G,2)$-arc-transitive\/} if $G\leq\Aut(\G)$ and $G_v$ is primitive or $2$-transitive  on $\G(v)$, respectively, for each $v\in V(\G)$. Similarly, when $G=\Aut(\G)$, a $G$-locally primitive or locally $(G,2)$-arc-transitive graph $\G$ is simply called {\em locally primitive\/} or {\em locally $2$-arc-transitive\/}, respectively.

A group $G$ of permutations of a set $V(\G)$ is called \emph{regular} if it is transitive, and some (and hence all) stabilisers $G_v=1$ are trivial.
More generally $G$ is called \emph{semiregular} if the stabiliser $G_v=1$ for all $v\in V(\G)$. So $G$ is regular if and only if it is semiregular and transitive.

\subsubsection{Normal quotients and normal covers of graphs}
Let $\G$ be a regular graph. Assume that $G\leq\Aut(\G)$ is such that $\G$ is $G$-vertex-transitive or $G$-edge-transitive. Let $N$ be a normal subgroup of $G$ such that $N$ is intransitive on $V(\G)$. The {\em $N$-normal quotient graph\/} of $\G$ is defined as the graph $\G_N$ with vertices the $N$-orbits in $V(\G)$ and with two distinct $N$-orbits adjacent if there exists an edge in $\G$ consisting of one vertex from each of these orbits. If $\G_N$ and $\G$ have the same valency, then we say that $\G$ is an {\em $N$-normal cover} of $\G_N$.

\subsubsection{Cayley graphs}\label{sub:cay}
Given a finite group $G$ and an inverse-closed subset $S\subseteq G\setminus\{1\}$ (that is, $s^{-1}\in S$ for all $s\in S$), the {\em Cayley graph} $\Cay(G,S)$ on $G$ with respect to $S$ is a graph with vertex set $G$ and edge set $\{\{g,sg\}\ :\ g\in G,s\in S\}$. For any $g\in G$ define
\[
R(g): x\mapsto xg\ \mbox{for $x\in G$ and set $R(G)=\{R(g)\ :\ g\in G\}$.}
\]
Then  $R(G)$ is a regular permutation group on $V(\G)$ (see, for example \cite[Lemma 3.7]{PS}) and is a subgroup of $\Aut(\Cay(G,S))$ (as $R(g)$ maps each edge $\{x,sx\}$ to an edge $\{xg,sxg\}$). For briefness, we shall identify $R(G)$ with $G$ in the following.
Let
\[
\Aut(G,S)=\{\a\in\Aut(G): S^\a=S\}.
\]
It was proved by Godsil~\cite{Godsil-1981} that the normaliser of $G$ in $\Aut(\Cay(G,S))$ is $G: \Aut(G,S)$. A Cayley graph $\Cay(G,S)$ is said to be {\em normal} if $G$ is normal in $\Aut(\Cay(G,S))$ (see \cite{Xu98}); this is equivalent to the condition $\Aut(\Cay(G,S))=G: \Aut(G,S)$. Following \cite{Praeger-NE-Cay}, we say that $\Cay(G,S)$ is {\em normal-edge-transitive} if $G:\Aut(G,S)$ is transitive on the edge set of $\Cay(G,S)$. Note that  $\Cay(G, S)$ is normal-edge-transitive if and only if either $\Aut(G,S)$ is transitive on $S$, or has two orbits in $S$ such that each is the set of  inverses of elements of the other (see \cite[Proposition~1(c)]{Praeger-NE-Cay}).

Cayley graphs are precisely those  graphs $\Gamma$ for which $\Aut(\G)$ has a subgroup that is regular on $V(\G)$. For this reason we say that a  graph $\Gamma$ is a \emph{bi-Cayley graph} if $\Aut(\G)$ has a subgroup $H$ which is semiregular on $V(\G)$ with two orbits.

\subsubsection{Notation for groups}
For a positive integer $n$,  $C_n$ denotes a cyclic group of order $n$, and $D_{2n}$ denotes a dihedral group of order $2n$.
For a group $G$, we denote by $1$, $Z(G)$, $\Phi(G)$, $G'$, $\soc(G)$ and $\Aut(G)$, the identity element, the centre, the Frattini subgroup, the derived subgroup, the socle and the automorphism group of $G$, respectively. For a subgroup $H$ of a group $G$,
denote by $C_G(H)$ the centralizer of $H$ in $G$ and by $N_G(H)$ the
normalizer of $H$ in $G$. For elements $a,b$ of $G$, the {\em commutator} of $a,b$ is defined as $[a,b]=a^{-1}b^{-1}ab$.
If $X,Y\subseteq G$, then $[X,Y]$ denotes the subgroup generated by all the commutators $[x,y]$ with $x\in X$ and $y\in Y$.

\subsection{Six useful lemmas}

The first result proves some properties of mixed dihedral groups (see Definition~\ref{mixed-dihedral-group}) mentioned in Section~\ref{s:intro}.

\begin{lem}\label{prop-of-mixed-dihedral}
Let $H$ be a mixed dihedral group relative to $X$ and $Y$, where $X, Y$ are subgroups of $H$ and $X\cong Y\cong C_2^n$. Then the following hold.
\begin{enumerate}
\item[{\rm (1)}] Let $\phi : H\to H/H'$ be the natural projection map given by $\phi : h\to hH'$. Then $H/H'\cong \phi(X)\times \phi(Y)\cong X\times Y$.

\item[{\rm (2)}] For all $h, g\in H$, $|Xh\cap Yg|\leq 1$.
 \end{enumerate}
\end{lem}

\f\demo (1)\ Let $X=\langle x_1,\ldots,x_n\rangle =C_2^n$ and $Y=\langle y_1,\cdots,y_n\rangle =C_2^n$. Then $H/H'=\phi(H)$ is generated by the $2n$ elements $\phi(x_i),\phi(y_i)$, for $1\leq i\leq n$, each of which has order at most $2$. Since by definition $|H/H'|=2^{2n}$, it follows that $\phi(X)\cong X$, $\phi(Y)\cong Y$, $H/H'\cong\phi(X)\times \phi(Y)$, and also $X\cap Y=1$.

(2)\ Suppose that $a, b\in Xh\cap Yg$. Then $Xa=Xh=Xb$, and so $ab^{-1}\in X$. Similarly, $Ya=Yg=Yb$, and so $ab^{-1}\in Y$. Thus, $ab^{-1}\in X\cap Y=1$, that is, $a=b$. Hence, $|Xh\cap Yg|\leq 1$.\hfill\qed

The second lemma concerns normal quotients of locally primitive graphs.

\begin{lem}\label{quot}
Let $\G$ be a connected regular $G$-locally primitive bipartite graph of valency $k>1$ and with bipartition $V(\G)=O_1\cup O_2$, so each $|O_i|>1$. Suppose that $N\unlhd G$ is such that $N$ fixes both $O_1$ and $O_2$ setwise, and $N$ is intransitive on $O_1$ and on $O_2$. Then
\begin{enumerate}
  \item [{\rm (1)}] $\G$ is an  $N$-normal cover of the quotient graph $\G_N$ of $\G$.
  \item [{\rm (2)}] $N$ acts semiregularly on $V(\G)$, $N$ is the kernel of the $G$-action on $V(\G_N)$, and $G/N\leq\Aut(\G_N)$.
  \item [{\rm (3)}] $\G_N$ is $G/N$-locally primitive. Furthermore, if $\G$ is locally $(G,2)$-arc-transitive, then $\G_N$ is locally $(G/N, 2)$-arc-transitive.

  \item [{\rm (4)}] For $N\leq H\leq G$, $H$ is regular on $E(\G)$  if and only if $H/N$ is regular on $E(\G_N)$.
\end{enumerate}
\end{lem}

\f\demo Parts (1)--(3) are proved in \cite[Lemma~5.1]{GLP-tamc}. 
Now we prove part (4). Note that, by part (2), $N$ is semiregular so each $N$-orbit in $V(\G)$ has size $|N|$, and by part (1), $\G$ is a cover of $\G_N$ so each edge $\{x^N, y^N\}$ of $\G_N$ corresponds to exactly $|N|$ edges of $\G$ and the subgroup $N$ acts regularly and faithfully on them. Suppose that $N\leq H\leq G$. Suppose first that $H$ is regular on $E(\G)$. Then $|H|=|E(\G)|=|N|\cdot|E(\G_N)|$, and $H$ is transitive on $E(\G_N)$. Hence the stabiliser $H_e$ in $H$ of an edge $e:=\{x^N, y^N\}$ of $\G_N$ has order $|H_e|=|H|/|E(\G_N)|=|N|$. Since $N\leq H$ and $N$ fixes the edge $e$, we have $N\leq H_e$, and it follows that $H_e=N$. Hence $H$ induces a regular group $H/N$ on $E(\G_N)$. Conversely suppose that $H/N$ is regular on $E(\G_N)$, so the stabiliser in $H$ of an edge $e:=\{x^N, y^N\}$ of $\G_N$ is equal to $N$. Since $N$ acts regularly and faithfully on the set of $|N|$ edges with end-points in $x^N, y^N$, it follows that $H$ acts regularly on $E(\G)$. \hfill\qed

The next result gives a basic property of bi-Cayley graphs.

\begin{lem}\label{bi-Cayley}
Let $\G$ be a connected bi-Cayley graph of a group $H$ such that neither of the two $H$-orbits in $V(\G)$ contains an edge of $\G$, and let $N=N_{\Aut(\G)}(H)$. Then, for each $v\in V(\Gamma)$, the stabiliser $N_v$ acts faithfully on $\Gamma(v)$.
\end{lem}

\f\demo
By the definition of a bi-Cayley graph, $H$ acts semiregularly on $V(\G)$ with two orbits,
say $U$ and $W$. We may assume that $U=\{h_0\mid h\in H\}$ and $W=\{h_1\mid h\in H\}$,
and that $H$ acts on $U$ and $W$ as follows:
\[
h_i^g=(hg)_i, \quad\mbox{for all  $h, g\in H$ and $i=0, 1$.}
\]
Since by assumption $U$ and $W$ contain no edges of $\G$, it follows that $\G$ is a regular bi-partite graph, and we may assume that
$E(\G)$ is the set of pairs $\{h_0,g_1\}$ with $h, g\in H$ such that $gh^{-1}$ lies in a certain subset $S$ of $H$, see $\cite[p.\, 505]{ZF}$.
In particular $\G(1_0)=\{s_1\mid s\in S\}\subseteq W$, and $\G(1_1)=\{(s^{-1})_0\mid s\in S\}\subseteq U$. Further, by \cite[Lemma~3.1]{ZF}, we may assume also that $1\in S$, and hence $1_1\in\G(1_0)$. Since $\G$ is connected it follows that  $H=\lg S\rg$.

Since $H$ has two orbits on $\Gamma$, in order to prove the lemma it is sufficient to show that $N_v$ is faithful on $\Gamma(v)$ for $v=1_0$ and $v=1_1$. By \cite[Theorem~1.1]{ZF}, we have $N_{1_0}=\{ \s_{\a,g}\ |\ \a\in\Aut(H), g\in H, S^\a=g^{-1}S\}$, where $\sigma_{\a,g}$ is defined as follows:
\begin{equation*}
\begin{aligned}
&~~\sigma_{\a,g}:~ h_{0}\mapsto(h^{\alpha})_{0},~ h_{1}\mapsto(gh^{\a})_{1},\quad\mbox{for all}\  h\in H,\\
\end{aligned}
\end{equation*}
and therefore $N_{1_01_1}=\{ \s_{\a,1}\mid \a\in\Aut(H), S^\a=S\}$.  Note that for  $\sigma_{\a,1}\in N_{1_01_1}$, we have, by definition, that $\sigma_{\a,1}:(s^{-1})_0\to (s^{-\alpha})_0$ and
$s_1\to (s^\alpha)_1$, for all $s\in S$.
Since $\G(1_0)=\{s_1\mid s\in S\}$ and $\G(1_1)=\{(s^{-1})_0\mid s\in S\}$, it follows that the kernels of the actions of $N_{1_0}$ on $\G(1_0)$, and of  $N_{1_1}$ on $\G(1_1)$, are equal to the same subgroup $K$, namely $K$ consists of all elements $\sigma_{\alpha,1}$ such that $\alpha$ fixes $S$ pointwise.  Since $H=\langle S\rangle$ this implies that $K$ is trivial. Thus the lemma is proved.
\hfill\qed

The last three lemmas are all related to complete bipartite graphs in some way. The first is a group theoretic characterisation.

\begin{lem}\label{lem:abelian-edge-regular}
 Let $\G$ be a connected graph with $|E(\G)|>1$, and suppose that $G\leq\Aut(\G)$ is abelian and edge-transitive. Then either
 \begin{enumerate}
     \item[{\rm (1)}] $\G$ is a cycle $\C_n$, for some $n\geq3$, and $G\cong C_n$ is vertex-transitive, or
     \item[{\rm (2)}] $\G=\K_{m,n}$ with biparts $\Delta, \Sigma$ of sizes $m,n$ respectively, and $G=M\times N<\Sym(\Delta)\times \Sym(\Sigma)$ with $|M|=m, |N|=n$ and $m+n>1$.
 \end{enumerate}
\end{lem}

\f\demo Since $|E(\G)|>1$ and $\G$ is connected, it follows that the group $G$ acts faithfully on $E(\Gamma)$. Choose an edge $\{u,v\}\in E(\G)$. Then since $G$ is edge-transitive and abelian, it follows that $G$ is regular on $E(\G)$ (see \cite[Theorem 3.2]{PS}), and so
$|G|=|E(\G)|$. Now each vertex $x$ lies in at least one  edge, say $\{x,y\}$ (since $\G$ is connected), and for some $g\in G$, $\{x,y\}=\{u,v\}^g=\{u^g,v^g\}$ (since $G$ is edge-transitive), and hence $x\in u^G$ or $x\in v^G$. Hence $G$ has at most two orbits in $V(\G)$.

Suppose first that $G$ is transitive on $V(\G)$. Then $\G$ is regular, say of valency $k$, and counting incident vertex-edge pairs we have $|V(\G)|\cdot k = |E(\G)|\cdot 2$.
Now $G$  is faithful on $V(\G)$ (by definition) and $G$ is abelian, and hence (again see \cite[Theorem 3.2]{PS})  $G$ is regular on $V(\G)$. Thus $|V(\G)|=|G|=|E(\G)$, and this implies that $k=2$. As $\G$ is connected, this means that $\G$ is a cycle $\C_n$ for some $n\geq3$, and the only abelian edge-transitive subgroup $G$ of $\Aut(\C_n)=D_{2n}$ is the cyclic group of rotations $C_n$, and part (1) holds.

Thus we may assume that $G$ has two vertex orbits, namely $\Delta=u^G$ of size $m$, and $\Sigma=v^G$ of size $n$. As $\G$ is connected, $\G$ is bipartite with biparts $\Delta$ and
$\Sigma$. Since $G$ is edge-transitive, it follows that $N:=G_u$ is transitive on $\G(u)$ (a subset of $\Sigma$), and as $G$ is abelian, the transitive $G$-action on $\Delta$ is regular so $N$ fixes $\Delta$ pointwise.
If $x\in\G_2(u)$, then $x\in\Delta$ and there exists $y\in\G(u)\cap \G(x)$. Then $G_x=N$ (as $G$ is regular on $\Delta$) and $G_x$ is transitive on $\G(x)$ (as $G$ is edge-transitive), and so $\G(x)$ is the $N$-orbit containing $y$, that is, $\G(x)=y^N=\G(u)$. This holds for every $x\in\G_2(u)$, and it follows, since $\G$ is connected that $\Sigma=\G(u)$ and that $N$ is transitive on $\Sigma$. The same argument with $v$ in place of $u$ proves that $\Delta=\G(v)$ and that $M:=G_v$ fixes $\Sigma$ pointwise and is transitive on $\Delta$. Thus $\G=\K_{m,n}$, and $M\times N\leq G$. Now $|G|=|E(\G)|=mn=|M\times N|$, and hence $G=M\times N$. Finally $m+n>1$ since  $|E(\G)|>1$.
%
\hfill\qed

We next record the result of a computer investigation of the two small graphs $\K_{4,4}$ and $\K_{8,8}$ using Magma~\cite{BCP}.

\begin{lem}\label{cor-val-16}
Let $n=2$ or $3$, and let $\G=\K_{2^n,2^n}$. If $G\leq \Aut(\G)$ is such that $\G$ is $G$-vertex-transitive and $G$-locally primitive, then $\G$ is $(G,2)$-arc-transitive  and $G$ contains a subgroup acting regularly on $V(\G)$.
\end{lem}

\begin{lem}\label{sufficient}
Let $\G$ be a connected $(G,2)$-arc-transitive graph, and let $u\in V(\G)$. Suppose that $\G$ is an $N$-normal cover of $\K_{2^n,2^n}$, for some normal $2$-subgroup $N$ of $G$. Then $\G$ is bipartite,  and one of the following holds:
\begin{enumerate}
  \item [{\rm (1)}] $\G$ is a normal Cayley graph of a $2$-group;
  \item [{\rm (2)}] $\G$ is a bi-Cayley graph of a $2$-group $H$ such that $G\leq N_{\Aut(\Gamma)}(H)$;
  \item [{\rm (3)}] $N\unlhd\Aut(\G)$.
\end{enumerate}
Moreover if the stabiliser $G_u$ acts non-faithfully on $\G(u)$, then part (3) holds.
\end{lem}

{
\f\demo
By assumption $\G_N\cong \K_{2^n,2^n}$, so $\G$ is bipartite, and also $\G$ is regular of valency $2^n$ since $\G$ covers $\G_N$. Now $\G$ is a $(G,2)$-arc-transitive graph with order a $2$-power (since $N$ is a $2$-group). We apply \cite[Theorem~1.1]{LMP2009} relative to the group $\Aut(\Gamma)$ and conclude that either $\G$ is a normal Cayley graph of a $2$-group, or $\Aut(\G)$ has a normal subgroup $M$ such that $\G_M\cong \K_{2^m,2^m}$ for some $m$ (see \cite[proof of Theorem~1.1 on p.~120]{LMP2009} and note that the graphs in case (ii) of that result do not arise as they are not bipartite).
In the former case, $G_u$ acts faithfully on $\G(u)$ (see, for example, {\cite[p.~4610]{L-TAMS-2006}}), and hence $\G$ satisfies part (1) and the lemma is proved in this case. Thus
we may assume that the latter holds. By Lemma~\ref{quot}~(1) and (2), $\G$ is an $M$-normal cover of $\G_M$, and $M$ is semiregular on $V(\G)$. In particular $\G$ has valency $2^m$ so that $m=n$, and $|V(\G)|=|V(\G_M)|\cdot|M|=2^{n+1}|M|$. As $\G$ is also an $N$-normal cover of $\K_{2^n,2^n}$, again by Lemma~\ref{quot}, $N$ is semiregular on $V(\G)$ and  $|V(\G)|=|V(\G_N)|\cdot|N|=2^{n+1}|N|$. It follows that $|N|=|M|$, and so $M$ is a $2$-group. Suppose first that $N=M$. Then part (3) holds and there is nothing further to prove.


Assume therefore that $N\neq M$, and recall that $N\lhd G$ and $M\lhd\Aut(\G)$. Since $\G_N\cong\G_M\cong \K_{2^n,2^n}$, both $M$ and $N$ stabilise each bipart of $\G$.
Let $G^+$ be the subgroup of $G$ stabilising each of the biparts of $\Gamma$. Then $NM\unlhd GM$, $NM\leq G^+M$ and $NM$ is a $2$-group.
Since $\G$ is $(G,2)$-arc-transitive, it follows from Lemma~\ref{quot}~(3) that $GM/M$ is $2$-arc-transitive on $\G_M$, and in particular,  $G^+M/M$ acts $2$-transitively on each part of size $2^n$ of $\Gamma_M=\K_{2^n,2^n}$.
Moreover, since $N\neq M$, $NM/M$ is a non-trivial normal $2$-subgroup of $G^+M/M$, and so $NM/M$ induces a regular action of $C_2^n$ on each part of  $\G_M$.
Suppose that $NM/M$ acts faithfully on each part of $\Gamma_M$. Then $N/(N\cap M)\cong NM/M\cong C_2^n$ (since each part of $\Gamma_M$ has size $2^n$), and $NM$ acts faithfully and regularly on each part of $\Gamma$ (since $M$ is semiregular on $V(\G)$). Thus $\Gamma$ is a  bi-Cayley graph of the $2$-group $H:=NM$, and as $H\unlhd GM$, $G$ is contained in the normaliser of $H$ in $\Aut(\Gamma)$, and part (2) holds. Also $G_u$ is faithful on $\G(u)$ by Lemma~\ref{bi-Cayley},  and there is nothing further to prove.

Therefore we may assume that  $NM/M$ is unfaithful on at least one of the parts of $\G_M$. Since $NM\unlhd GM$ and $GM$ is transitive on $V(\G)$, the action of $NM$ is unfaithful on each of the parts of $\G_M$. We will derive a contradiction, and this will complete the proof (even of the last assertion) of the lemma.

We showed above that the group induced by $NM/M$ on each of the parts of $\G_M$ is $C_2^n$,  and hence $NM/M$ is isomorphic to a subdirect subgroup of $C_2^n\times C_2^n$. In particular $NM/M$ is elementary abelian.
Let $U,W$ be the two parts of the bipartition for $\G_M$, so the actions of $NM/M$ on $U$ and on $W$ are not faithful. We may assume that the vertex $u$ lies in a part $\bf u$ of $U$.  Let $K/M$ be the kernel of $NM/M$ acting on $U$, and note that $K/M$ acts faithfully on $W$. Since $NM/M$ is abelian and regular on $U$, $K/M$ is the stabiliser of ${\bf u}$ in $NM/M$. So $K/M=(NM/M)_{\bf u}\unlhd (GM/M)_{\bf u}$. Since $GM/M$ is $2$-arc-transitive on $\G_M=\K_{2^n,2^n}$, $(GM/M)_{\bf u}$ is $2$-transitive on the part $W$ of $\G_M$.
It follows that its nontrivial normal $2$-subgroup $K/M=(NM/M)_{\bf u}$ is abelian and regular on $W$, and we conclude that $NM/M\cong C_2^n\times C_2^n=C_2^{2n}$, and that $NM/M$ acts regularly on the edge set of $\G_M$. It follows from Lemma~\ref{quot}~(4) that $NM$ is regular on the edges of $\G$.

Now $NM/M\cong C_2^{2n}$ implies that the Frattini subgroup $\Phi(NM)\leq M$, and since $NM\unlhd GM$, also $\Phi(NM)\unlhd GM$. Since $GM$ is $2$-arc-transitive on $\G$, by Lemma~\ref{quot}~(1)--(3), $\G$ is a $\Phi(NM)$-normal cover of the quotient graph $\G_{\Phi(NM)}$ with $GM/\Phi(NM)$ as a $2$-arc-transitive group of automorphisms. In particular $\G_{\Phi(NM)}$ has valency $2^n$. Since $NM$ is regular on the edge set of $\G$, by Lemma~\ref{quot}~(4), $NM/\Phi(NM)$ is also regular on the edge set of $\G_{\Phi(NM)}$, and since $NM/\Phi(NM)$ is elementary abelian {(see \cite[5.2.12]{Robinson})}, it follows from Lemma~\ref{lem:abelian-edge-regular} that $\G_{\Phi(NM)}\cong \K_{2^n,2^n}$, and so $NM/\Phi(NM)\cong C_2^{2n}$. This implies that $M=\Phi(NM)$.  Since $NM$ is generated by $N\cup M=N\cup\Phi(NM)$ it follows that $NM$ is generated by $N$ ( see \cite[5.2.12]{Robinson}), and hence $NM=N$.  Thus $M\leq N$. However we proved above that $|M|=|N|$, so we conclude that $M=N$, which contradicts $NM/M\cong C_2^{2n}$.\hfill\qed
}


\section{Proof of Theorem~\ref{is-affine}}

The goal of this section is to prove Theorem~\ref{is-affine}.
The proof will depend heavily on the following result, Proposition~\ref{basic-auto}. In turn, the proof of  Proposition~\ref{basic-auto} relies on the classification of the finite $2$-transitive permutation groups of $2$-power degree, and hence on the finite simple group classification.

\begin{prop}
\label{basic-auto}
Let $n\geq 2$ be a positive integer and let $\G=\K_{2^n,2^n}$ with biparts $U$ and $W$. Let $G\leq\Aut(\G)$ be $2$-arc-transitive on $\G$, and let $G^+=G_U=G_W$, the setwise stabilizer in $G$ of each of the biparts. Then one of the following holds.
\begin{enumerate}
  \item [{\rm(1)}]\ $\soc(G^+)\cong A_{2^n}\times A_{2^n}$ with $n\geq3$;
  \item [{\rm(2)}]\ $\soc(G^+)\cong\PSL(2,p)\times\PSL(2,p)$ with $p=2^n-1$ a prime and $n\geq 3$;
  \item [{\rm(3)}]\ $\soc(G^+)\cong C_2^{2n}$,  acting regularly on the edge set of $\G$;
  \item [{\rm(4)}]\ $n=3$ and $G^+={\rm AGL(3,2)}, G_u=\PSL(2,7), G_w=\SL(3,2)$, with $G_{u}\cap G_w=\mz_7:\mz_3$, where $u\in U$ and $w\in W$.
\end{enumerate}
Moreover, if $G$ does not contain a  subgroup acting regularly on $V(\G)$, then $n\geq 4$ and
case {\rm(3)} holds.
\end{prop}



{
\demo
Let $u\in U$ and $w\in W$. Since $G$ is $2$-arc-transitive on $\G$, there exists an element $g\in G$ such that $g$ interchanges $u$ and $w$, and replacing $g$ by an odd power of $g$ if necessary, we may assume that $g$ has order a $2$-power. Then $g^2\in G_{u}\cap G_w\subseteq G^+$. Also, since $\G$ is $(G,2)$-arc-transitive, the subgroup $G_u$ is $2$-transitive on $\G(u)=W$ and hence also $G^+$ is $2$-transitive on $W$. Similarly $G^+$ is 2-transitive on $U=\G(w)$. Thus $G^+$ is a subdirect subgroup of $H\times H^g$ for some $2$-transitive subgroup $H$ of $\Sym(U)=S_{2^n}$, and $G_w^U, G_u^W$ are $2$-transitive subgroups of $H, H^g$, respectively.

Suppose first that $n=2$. Then $H=A_4$ or $S_4$, and since $G_u^W$ is $2$-transitive, $\soc(G_u^W)=C_2^2$. Hence $\soc(G^+)=(C_2^2)^2=C_2^{2n}$ as in part (3). Moreover, it follows from Lemma~\ref{cor-val-16} that $G$ contains a subgroup acting regularly on $V(\G)$. Thus we may assume that $n\geq3$.  Next suppose that case (4) holds for $G$. We note that a subgroup $G$ of $\Aut(\G)$ with all the required properties exists, see~\cite[Theorem (b), and pp.~26--27]{praeger1982} which presents a construction due to Chris Rowley.  In particular, by Lemma~\ref{cor-val-16}, $G$ contains a subgroup acting regularly on $V(\G)$.  Thus we may assume further that $G$ is not as in case (4).
It now follows from {\rm\cite[Corollary~1.2]{FanLLP2013}} that
\begin{center}
    $G^+=(M\times M^g).P$, where $M$ fixes $W$ pointwise\\ and $M.P$ is faithful and $2$-transitive on $U$.
\end{center}
Let $T$ be the socle of the $2$-transitive group $M.P$ so that $T\leq M<\Sym(U)$ and either $T$ is elementary abelian or $T$ is a nonabelian simple group (by a  theorem of Burnside, see \cite[Theorem 3.21]{PS}), and as the degree $|U|=2^n$ it follows from \cite[Theorem~2.2]{LMP2009} that either $T\cong C_2^n$ with $M.P\lesssim {\rm AGL}(n,2)$, or $T\cong A_{2^n}$ (with $n\geq3$), or $T\cong \PSL(2,p)$ with $p=2^n-1$ a prime and $n\geq 3$.
In each of these cases $T$ is a minimal normal subgroup of $G^+$ which is transitive on $U$ and fixes $W$ pointwise, and also $T^g$ is a minimal normal subgroup of $G^+$ with $T^g$  transitive on $W$ and fixing $U$ pointwise. Thus $T\times T^g$ is contained in $\soc(G^+)$. If $\soc(G^+)$ is strictly larger than $T\times T^g$ then $G^+$ has a minimal normal subgroup $N$ such that $N\cap (T\times T^g)=1$, and in this case $N\leq C_{G^+}(T\times T^g)$.
However $C_{G^+}(T\times T^g)\leq C_{\Sym(U)}(T)\times C_{\Sym(W)}(T^g)$. If $T=C_2^n$ then $T$ is self-centralising in $\Sym(U)$ by \cite[Theorem 3.6]{PS}, while if $T\cong A_{2^n}$, or $T\cong \PSL(2,p)$ then $T$ is $2$-transitive on $U$ and so  $C_{\Sym(U)}(T)=1$ by \cite[Theorem 3.2]{PS}. In either case it follows that $N\leq C_{G^+}(T\times T^g)\leq T\times T^g$. This is a contradiction, and hence $\soc(G^+)=T\times T^g$. Thus one of the cases (1), (2) or (3) holds.

To complete proof of the last assertion, we need to prove that $G$ has a subgroup $X$ which is regular on $V(\G)$ if case (1) or (2) holds. This follows from Lemma~\ref{cor-val-16} if $n=3$, so we may assume that $n\geq4$.
Therefore we have $\soc(G^+)=T\times T^g$, where either $T\cong A_{2^n}$, or $\PSL(2,p)$ with $p=2^n-1$ a prime.
Suppose first that there exists an involution in $G\setminus G^+$. Then we may assume that $|g|=2$.
In either case,  $T$ has a subgroup $H$ acting regularly on $U$; for example, $H=C_2^n$ or $H=D_{p+1}$ respectively.
Then $H^g\leq T^g$ and $H^g$ acts regularly on $W$ and fixes $U$ pointwise. Moreover $H\times H^g\leq G^+$ and so $L:=\{hh^g\mid h\in H\}$ is a subgroup of $G^+$ that is semiregular on $V(\G)$ with orbits $U$ and $W$. Since, for each $h\in H$,  we have  $(hh^g)^g= h^gh= hh^g$, it follows that $L\times\lg g\rg$ is regular on $V(\G)$. Thus  we may assume that
$G\setminus G^+$ contains no involutions.

Now $\soc(G^+)=T\times T\unlhd G\leq \mathcal{W} =H\wr S_2$, with the pair $(T,H)=(A_{2^n},S_{2^n})$ or $(\PSL(2,p), \PGL(2,p))$, and in both cases $\mathcal{W}/\soc(G^+)\cong D_8$ and $G/\soc(G^+)\ne 1$. Further, in both cases there exists an involution $y\in H\setminus T$ (with $H$ acting on $U$ and fixing $W$ pointwise) and an involution $z\in \mathcal{W}\setminus(H\times H)$ (generating the top group of the wreath product) such that $\langle y, z\rangle\cong D_8$. Hence $G=\soc(G^+)\rtimes K$ for some nontrivial subgroup $K$ of  $\langle y, z\rangle$ with $K\not\leq \langle y\rangle\times\langle y^z\rangle$. The condition that  $G\setminus G^+$ contains no involutions implies that $K=\langle yz\rangle\cong C_4$.

We now construct, in both cases, explicit subgroups of $G$ that are regular on $V(\G)$.
Let us take $U=\{1,2,\ldots, 2^n\}$ and $W=\{1', 2', \ldots, (2^n)'\}$, so that the generator $z$ of the top group of $\mathcal{W}$ is
\[
z=(1,1')(2,2')\cdots(2^n, (2^n)').
\]
Both $H=S_{2^n}$ and $H=\PGL(2,p)$ contain a $2^n$-cycle $x$ acting regularly on $U$, and an involution $y$ that inverts $x$. Hence,  replacing  $x, y$ (and hence also $H$ in the second case) by conjugates in $S_{2^n}$ if necessary, we may assume that
\[
x=(1,2,\ldots, 2^n)\ \mbox{and}\ y=(2, 2^n)(3, 2^n-1)\cdots(2^{n-1}, 2^{n-1}+2).
\]
We note that $x, y\in H\setminus T$ so the product $xy\in T$, and $xy$ is an involution with $\langle x, y\rangle \cong D_{2^{n+1}}$. Further, taking $x, y$ to act on $U$ and fix $W$ pointwise, we see that $y, z$ are as in the previous paragraph, that is to say, $\lg y, z\rg\cong D_8$ and $\mathcal{W}=\soc(G^+)\rtimes \langle y,z\rangle$. Thus, as argued above, $G=\soc(G^+)\rtimes K$ with $K=\langle yz\rangle\cong C_4$. Finally we note that $xz = (xy)(yz)$ lies in $G$ since  $xy\in T\times 1\leq \soc(G^+)$ and $yz\in K\leq G$. Computing the product $xz$ explicitly we find that
\[
xz= (1,2,\ldots, 2^n) \cdot (1,1')(2,2')\cdots (2^n, (2^n)') = (1',1,2',2,\dots (2^n)',2^n).
\]
Thus $\langle xz\rangle \cong C_{2^{n+1}}$ and acts regularly on $V(\G)$. This completes the proof.
\hfill\qed
}

\medskip
\f{\bf Proof of Theorem~\ref{is-affine}}\ Let $\G, G, N$ be as in the statement of Theorem~\ref{is-affine}. Then, by \cite[Theorem 4.1]{ONS}, $G/N$ is $2$-arc-transitive on $\G_N=\K_{2^n,2^n}$. Let $G^+/N$ be the subgroup of $G/N$ stabilising both parts of the bipartition of $\G_N=\K_{2^n,2^n}$. Suppose that $\G$ is not a Cayley graph. Then $G$ does not contain a subgroup acting regularly on $V(\G)$, and it follows from the last assertion of Proposition~\ref{basic-auto} that $n\geq4$, and part (3) of Proposition~\ref{basic-auto} holds, that is, $\soc(G^+/N)\cong C_2^{2n}$, acting regularly on the edge set of $\G_N$. In other words,
$\G_N$ is $G/N$-edge-affine. \hfill\qed

\section{Characterisation}

The goal of this section is to prove  Theorem~\ref{iff-characterization}.

Throughout this section, we let $H$ be an $n$-dimensional mixed dihedral group relative to $X$ and $Y$ with $|X|=|Y|=2^n\geq4$, and let $C(H,X,Y)$ and $\Sigma(H,X,Y)$ be the graphs  defined in Definition~\ref{mixed-dihedrant}.
We say that an edge $\{h,g\}$ of $C(H,X,Y)$ is an \emph{$X$-edge} if $hg^{-1}\in X$, or a \emph{$Y$-edge} if $hg^{-1}\in Y$. By Lemma~\ref{prop-of-mixed-dihedral}~(2),  $X\cap Y=1$. Hence these concepts are well defined and each edge of $C(H,X,Y)$ is either an $X$-edge or a $Y$-edge. Note that, by definition, each $X$-edge is of the form $\{g,xg\}$ and each $Y$-edge is $\{g,yg\}$, for some $g\in H, x\in X\setminus\{1\}$, and $y\in Y\setminus\{1\}$. In our first result we describe several graph theoretic links between $C(H,X,Y)$ and $\Sigma(H,X,Y)$.


\begin{lem}\label{lem:prop-mixed-dih}
Let $H, X, Y, n$ be as above. Then the following hold.
\begin{enumerate}
\item [{\rm (1)}] For each triangle $(3$-clique$)$ $\{g,h,k\}$ in $C(H,X,Y)$,  either all three edges are $X$-edges or all three edges are $Y$-edges.

\item [{\rm (2)}] Each $X$-edge $\{g, xg\}$ of  $C(H,X,Y)$ lies in a unique maximal clique, namely $Xg$, and each $Y$-edge $\{g,yg\}$ lies in a unique maximal clique, namely $Yg$.

\item [{\rm (3)}] $\Sigma(H,X,Y)$ is the clique graph of $C(H,X,Y)$.

\item [{\rm (4)}] The map $\varphi:z\to \{Xz,Yz\}$, for $z\in H$, is a bijection $\varphi:H\to E(\Sigma(H,X,Y))$, and induces a graph isomorphism from  $C(H,X,Y)$ to the line graph $\mathcal{L}(\Sigma(H,X,Y))$ of $\Sigma(H,X, Y)$.

\item [{\rm (5)}]  Moreover, $\Aut(C(H,X,Y))=\Aut(\Sigma(H,X,Y))= \Aut(\mathcal{L}(\Sigma(H,X,Y)))$.
\end{enumerate}
\end{lem}

\f\demo (1)\quad For convenience, we let $\Gamma:= C(H,X,Y)$ and $\Sigma:=\Sigma(H,X,Y)$. Suppose that $\{g,h\}$ is an $X$-edge so $x:=hg^{-1}\in X\setminus\{1\}$, and that $\{h,k\}$ is a  $Y$-edge so $y:=kh^{-1}\in Y\setminus\{1\}$. Then $kg^{-1} =yx\not\in (X\cup Y)\setminus\{1\}$, so $\{g,k\}$ is not an edge.  This implies part (1).

(2)\quad It follows from the definition of $\G$ in Definition~\ref{mixed-dihedrant} that each pair $\{ xg, x'g\}$ of distinct elements of $Xg$ is an $X$-edge, since $(x'g)(xg)^{-1}=x'x^{-1}\in X\setminus\{1\}$. Hence $Xg$ is a clique and it contains $\{g, xg\}$. Also, by part (1),  each clique of size at least $3$ containing $\{g, xg\}$ must contain only $X$-edges. The only $X$-edges incident with $g$ or $xg$ are of the form $\{g, x'g\}$ or $\{xg, x'xg\}$, respectively, for some $x'\in X\setminus\{1\}$, and hence each such clique is contained in $Xg$. Thus, $Xg$ is a maximal clique of $\G$, and is the unique maximal clique containing $\{g, xg\}$. Similarly, $Yg$ is a maximal clique and is the unique maximal clique containing the $Y$-edge $\{g,yg\}$. This proves part (2).

(3)\quad By part (2), the clique graph $\Sigma(\G)$ of $\G$ has  vertex set $\{Xg\mid g\in H\}\cup \{Yg\mid g\in H\}$, and this is equal to $V(\Sigma)$ by Definition~\ref{mixed-dihedrant}. Moreover, two maximal cliques are adjacent in $\Sigma(\G)$ if and only if they contain at least one common vertex. Since distinct cosets of $X$, or of $Y$ are disjoint, each edge in the clique graph $\Sigma(\G)$ is of the form $\{Xg, Yh\}$ such that $g,h\in H$ and $Xg\cap Yh\neq\emptyset$.  Thus, the edges of the clique graph $\Sigma(\G)$ are precisely the edges of $\Sigma$, and so $\Sigma=\Sigma(\G)$, proving part (3).

(4)\quad The fact that $\Gamma$ is isomorphic to the line graph of $\Sigma$ can be deduced from \cite[Corollary 1.6]{DJLP-JCTA}, since the  subgraph of $\G$ induced on the neighbourhood $\Gamma(1)=(X\cup Y)\setminus\{1\}$ is isomorphic to $2\K_{2^n-1}$, where $|X|=|Y|=2^n$. However we require an explicit isomorphism for our later work.

By Definition~\ref{mixed-dihedrant}, for each $z\in H$, the image $\varphi(z)$ is an edge of $\Sigma$ since $\varphi(z)=\{Xz, Yz\}$ and $Xz\cap Yz = (X\cap Y)z = \{z\}\ne\emptyset$. Thus $\varphi:H\to E(\Sigma)$ is well defined. Next, if $z,w\in H$ and $\varphi(z)=\varphi(w)$, then $Xz=Xw$ and $Yz=Yw$, and hence $zw^{-1}\in X\cap Y=\{1\}$, so $z=w$. Thus $\varphi$ is one-to-one. Next, it follows from Lemma~\ref{prop-of-mixed-dihedral}(2) that, for each edge $e=\{Xg, Yh\}$ of $\Sigma$ the defining property of an edge, namely $Xg\cap Yh\ne\emptyset$, is equivalent to  $Xg\cap Yh=\{z\}$ for some unique $z\in H$. Hence the edge can be written as $e=\{Xz, Yz\}$. It follows that $\varphi$ is onto, and hence is a bijection.

Now $\varphi$ induces a natural bijection from the set of unordered pairs from $H$ to the set of unordered pairs of edges of $\Sigma$, namely $\varphi:\{h,g\}\to \{\varphi(g),\varphi(h)\}$. For an $X$-edge $\{g,xg\}$, the images $\varphi(g), \varphi(xg)$ share a common $X$-coset, namely $Xg=Xxg$. Similarly, for a $Y$-edge $\{h,yh\}$, the images  $\varphi(h), \varphi(yh)$ share a common $Y$-coset, namely $Yh=Yyh$. Thus, for arbitrary $g,h\in H$, each of $\varphi(\{g,xg\})$ and $\varphi(\{h,yh\})$ is an edge-pair from $\Sigma$ (hence a vertex-pair from $\mathcal{L}(\Sigma)$) which intersects nontrivially, and hence forms an edge of $\mathcal{L}(\Sigma)$. Thus the restriction of this induced map $\varphi$ to $E(\Gamma)$ is a one-to-one map into $E(\mathcal{L}(\Sigma))$. We claim that this restriction is onto.
Since $\varphi$ on edges is one-to-one, we have $|\varphi(E(\G))|=|E(\G)|=|V(\G)|\cdot|\G(1)|/2 = |H|\cdot(2^n-1)$. Also, since $\varphi$ (acting on $H$) is  bijective, $|V(\mathcal{L}(\Sigma))| = |E(\Sigma)| = |\varphi(H)|= |H|$. Now by Remark~\ref{remark-sigma-graph}, $\Sigma$ is regular of valency $2^n$, and hence $\mathcal{L}(\Sigma)$ is regular of valency $2(2^n-1)$. It follows that
$|E(\mathcal{L}(\Sigma))| = |H|\cdot2(2^n-1)/2 = |\varphi(E(\G))|$. Thus the induced map $\varphi$ on edges is onto, proving the claim. Therefore $\varphi$ induces a bijection on the edge sets, and hence $\varphi$ is a graph isomorphism from $\G$ to the line graph $\mathcal{L}(\Sigma)$, completing the proof of part (4).

(5)\quad Considering the induced action of $\Aut(\G)$ on the set of maximal cliques of $\G$, and noting that adjacent cliques $Xg$, $Yh$ intersect in a single vertex of $\Gamma$ (by Lemma~\ref{prop-of-mixed-dihedral}) it follows $\Aut(\G)$ induces a faithful action as a subgroup of automorphisms of the clique graph of $\G$. Hence by part (3), $\Aut(\G)\leq \Aut(\Sigma)$. Also, considering the  induced action of $\Aut(\Sigma)$ on the set of edges of $\Sigma$, and noting that adjacent edges in the line graph of $\Sigma$ intersect in a unique vertex of $\Sigma$, it follows $\Aut(\Sigma)$ induces a faithful action as a subgroup of automorphisms of $\mathcal{L}(\Sigma)$. Hence  $\Aut(\Sigma)\leq \Aut(\mathcal{L}(\Sigma))$, so we have
$\Aut(\G)\leq \Aut(\Sigma)\leq \Aut(\mathcal{L}(\Sigma))$.
It follows from part (4) that equality holds, proving part (5).
 \hfill\qed


Now we consider the symmetry of these graphs in more detail. Recall that $H'$ denotes the derived subgroup of $H$.

\begin{lem}\label{lem:mixed-dih2}
Let $H, X, Y, n$ be as above, let $\Sigma =\Sigma(H,X,Y)$, and let $G=H:A(H,X,Y)$, where $A(H,X,Y)$ is the setwise stabiliser in $\Aut(H)$ of $X\cup Y$. Then the following hold.
\begin{enumerate}
\item [{\rm (1)}] The group $G$ acts as a subgroup of automorphisms on $\Sigma$ as follows, for $h,z\in H, \sigma\in A(H,X,Y)$, and $\varphi:H\to E(\Sigma)$ as in Lemma~\ref{lem:prop-mixed-dih}(4):
\begin{align*}
    \mbox{Vertex action:} && h:Xz\to Xzh,\ &  Yz\to Yzh && \sigma:Xz\to X^\sigma z^\sigma,\ \  Yz\to Y^\sigma z^\sigma\\
    \mbox{Edge action:} &&  h:\varphi(z)\to \varphi(zh) &&& \sigma:\varphi(z)\to \varphi(z^\sigma)\\
\end{align*}
The subgroup $H$ acts regularly on $E(\Sigma)$ and has two orbits on $V(\Sigma)$. In particular, this $G$-action is edge-transitive.

\item [{\rm (2)}] The $H'$-normal quotient graph $\Sigma_{H'}$ of $\Sigma$ is isomorphic to $\K_{2^n,2^n}$ and admits $G/H'$ acting faithfully as an edge-transitive group of automorphisms. Moreover, $\Sigma$ is an $H'$-normal cover of $\K_{2^n,2^n}$.

\item [{\rm (3)}] $A(H, X, Y)\leq (\Aut(X)\times\Aut(Y)): C_2\cong(\GL(n,2)\times\GL(n,2)): C_2$.

\item [{\rm (4)}] Let $H\unlhd L\leq G$. If $C(H,X,Y)$ is $L$-edge-transitive, then $\Sigma(H,X,Y)$ is $(L,2)$-arc-transitive.
\end{enumerate}
\end{lem}

\f\demo
(1) \quad By Lemma~\ref{lem:prop-mixed-dih}(3), $\Sigma$ is the clique graph of $\Gamma:=C(H,X,Y)$ and hence $\Aut(\G)$ induces a subgroup of automorphisms of $\Sigma$ via its induced action on subsets of $V(\G)$ (namely on the maximal cliques of $\G$). In particular, since  by definition $\G$ is a Cayley graph for $H$ we have $G\leq \Aut(\G)$, and this gives a $G$-action as a subgroup of automorphisms of $\Sigma$.  Since $H$ acts by right multiplication, and $A(H,X,Y)$ acts naturally, on $V(\G)=H$, the vertex-actions of $h, \sigma$ are as in the statement, and in particular $H$ has two vertex-orbits, namely $[H:X]=\{ Xz\mid z\in H\}$ and $[H:Y]=\{ Yz : z\in H\}$. (Note that $\sigma$ fixes the set $\{X,Y\}$ setwise.) Also, by Lemma~\ref{lem:prop-mixed-dih}(4), $\varphi:z\to \{Xz, Yz\}$ is a bijection $H\to E(\Sigma)$, and it follows from the definition of the $G$-action that an edge $\{Xz, Yz\}$ is mapped by $h, \sigma$ to $\{Xzh, Yzh\}$, $\{X^\sigma z^\sigma, Y^\sigma z^\sigma\}$, respectively, and hence $\varphi(z)^h=\varphi(zh)$ and $\varphi(z)^\sigma = \varphi(z^\sigma)$. Thus the edge-action is as asserted. In particular, as $H$ acts regularly by right multiplication on $H$, it follows that $H$ acts regularly on $E(\Sigma)$, and hence also $G$ is transitive on $E(\Sigma)$. This proves part (1).

(2)\quad Since $H'$ is normal in $H$, elements of $H$ permute the $H'$-orbits in $V(\Sigma)$ by right multiplication. Let $\Delta^X=\{Xh : h\in H'\}$ and $\Delta^Y=\{Yh : h\in H'\}$. Then $\Delta^X, \Delta^Y$ are the $H'$-orbits in $V(\Sigma)$ containing $X$ and $Y$, and they lie in the $H$-orbits $[H:X], [H:Y]$ on vertices, respectively. Further, it follows from Lemma~\ref{prop-of-mixed-dihedral}~(1) that each $g\in H$ is of the form $g=zxy$ for some $z\in H', x\in X, y\in Y$.

We claim that $\Delta^Xg=\Delta^Xy$, and that the number of $H'$-orbits in $[H:X]$ is $2^n$.
By definition, $\Delta^Xz=\Delta^X$. Also for each $Xh\in\Delta^X$ (where $h\in H'$), $Xhx = X(xhxh^{-1})h =X[x,h^{-1}]h$, and since $[x,h^{-1}]\in H'$, it follows that $Xhx=Xz'$ with $z'=[x,h^{-1}]h\in H'$, and hence $Xhx\in \Delta^X$. Thus also $\Delta^Xx=\Delta^X$, and so $\Delta^Xg=\Delta^Xy$.
Therefore each of the $H'$-orbits in $[H:X]$ is of the form $\Delta^Xy$ for some $y\in Y$.
Since $|Y|=2^n$, in order to prove that the number of $H'$-orbits in $[H:X]$ is $2^n$, it is sufficient to prove that $\Delta^Xy\ne \Delta^Xy'$ for distinct $y, y'\in Y$, or equivalently, to prove that $\Delta^Xy=\Delta^X$ (with $y\in Y$) implies that $y=1$. So suppose that $\Delta^Xy=\Delta^X$ for some $y\in Y$. Then $Xy\in \Delta^X$, so $Xy=Xz$ for some $z\in H'$, and hence  $z=xy$ for some $x\in X$. Using the map $\phi$ from Lemma~\ref{prop-of-mixed-dihedral}~(1), this implies that $\phi(x)=\phi(zy^{-1}) = \phi(y^{-1})\in\phi(X)\cap \phi(Y)=1$, and hence that $x=y=1$. Thus the claim is proved.

An analogous proof shows that $H'$ has exactly $2^n$ orbits in $[H:Y]$, and that these orbits are $\Delta^Yx$ for $x\in X$. Thus, the  $H'$-normal quotient of $\Sigma$ has $2\times 2^n$ vertices and, for each $x\in X$, by the definition of $\Sigma$, the orbit $\Delta^X$ is adjacent to $\Delta^Yx$  since by Lemma~\ref{lem:prop-mixed-dih}~(1) the cliques  $X$ and $Yx$ both contain $x$. Similarly  $\Delta^Y$ is adjacent to $\Delta^Xy$ for all $y\in Y$. Since the right multiplication action of $H$ induces a subgroup of automorphisms of the $H'$-normal quotient graph $\Sigma_{H'}$, and in this action $H$ is transitive on each of the two sets $[H:X]$ and $[H:Y]$, it follows that  $\Sigma_{H'}$ is isomorphic to $\K_{2^n,2^n}$.
Since $\Sigma$  and  $\K_{2^n,2^n}$ both have valency $2^n$, it follows that $\Sigma$ is an $H'$-normal cover of $\K_{2^n,2^n}$. Further, since $\Sigma$ is connected and a cover of $\Sigma_{H'}\cong\K_{2^n,2^n}$,  the kernel of the $G$-action on $V(\Sigma_{H'})$ is semiregular on $V(\Sigma)$,  and hence is equal to $H'$. Thus $G/H'$ acts faithfully as an edge-transitive group of automorphisms of $\Sigma_{H'}$. This proves part (2).

(3) \quad Now $G\leq\Aut(\Gamma)$, and its subgroup $A(H,X,Y)$ is the stabiliser in $G$ of the vertex $1$ of $\Gamma$. Since $X\cup Y$ generates $H$, $A(H,X,Y)$ acts faithfully on $\Gamma(1)=(X\cup Y)\setminus\{1\}$. It then follows from Lemma~\ref{lem:prop-mixed-dih}~(2) that $A(H,X,Y)$ acts faithfully on the set of maximal cliques of $\Gamma$, and from Lemma~\ref{lem:prop-mixed-dih} parts~(1) and~(2) that each element of $A(H,X,Y)$ either fixes setwise or interchanges the subsets $X\setminus\{1\}$ and $Y\setminus\{1\}$ of $\Gamma(1)$. The subgroup, of  $A(H,X,Y)$, of index 1 or 2,  fixing setwise each of these subsets, induces automorphisms of $X$ and $Y$, and part (3) follows.

(4)\quad Suppose that $H\leq L\leq G$ such that $L$ is transitive on $E(\G)$, and let $S=(X\cup Y)\setminus\{1\}$, and $s\in S$. Since $H$ is transitive on $V(\G)$, and $L$ is transitive on $E(\G)$, it follows from \cite[Proposition 1]{Praeger-NE-Cay} that the arc set of $\G$ is $(1,s)^L\cup (s,1)^L$. Moreover, since $s$ is an involution, right multiplication $R(s)\in H\leq L$ by $s$ maps $(1,s)$ to $(s,s^2)=(s,1)$, and hence $(1,s)^L=(s,1)^L$. Thus $\G$ is $L$-arc-transitive. In particular, the stabiliser $L_1$ of $1\in H$ acts transitively on $S$, and since $L_1\leq A(H,X,Y)$, it follows from part (3) that $L_1$ acts imprimitively on $S$ and each element of $L_1$ either fixes setwise each of $X\setminus\{1\}$ and $Y\setminus\{1\}$ or interchanges these two sets. Let $L_1^+$ denote the index 2 subgroup of $L_1$ fixing each of $X$ and $Y$ setwise. Then $L_1^+$ is transitive on each of the sets $X\setminus\{1\}$ and $Y\setminus\{1\}$.

Now $L$ is transitive on $V(\Sigma)$ by part (1) (since some element of $L_1$ interchanges $X\setminus\{1\}$ and $Y\setminus\{1\}$). The stabiliser in $L$ of the clique $X$ (a vertex of $\Sigma$), contains the subgroup $X$, which acts transitively on the set of cliques $Yx$ ($x\in X$) adjacent to $X$. Hence $\Sigma$ is $L$-arc-transitive. Also the stabiliser in $L$ of the arc $(X,Y)$ of $\Sigma$ contains $L_1^+$, which acts transitively on the set of $2$-arcs $(X, Y, Xy)$ ($y\in Y\setminus\{1\}$) extending this arc.  It follows that $\Sigma$ is $(L,2)$-arc-transitive, proving part (4).
\hfill\qed


\f{\bf Proof of Theorem~\ref{iff-characterization}~}
Suppose that $n$ is an integer with $n\geq2$, that $\Sigma$ is a graph, and $N\unlhd G\leq\Aut(\Sigma)$.
First, we prove that part (b) implies part~(a). Suppose, as in part (b), that $G$ has a normal subgroup $H$, where $H$ is an $n$-dimensional mixed dihedral group relative to $X$ and $Y$ with derived subgroup $H'=N$, and that the line graph $\mathcal{L}(\Sigma)$ of $\Sigma$ is $C(H,X,Y)$, and $C(H,X,Y)$ is $G$-edge-transitive. Lemma~\ref{lem:prop-mixed-dih}~(4) gives an explicit isomorphism $\varphi$ from
$C(H,X,Y)=\mathcal{L}(\Sigma)$ to the line graph $\mathcal{L}({\hat \Sigma})$ of ${\hat \Sigma}=\Sigma(H,X,Y)$. Since the vertex sets of $\mathcal{L}(\Sigma)$ and $\mathcal{L}({\hat \Sigma})$ are $E(\Sigma)$ and $E({\hat \Sigma})$, respectively, $\varphi$ is a bijection $E(\Sigma)\to E({\hat \Sigma})$, and similarly $\varphi$ yields a bijection $V(\Sigma)\to V({\hat \Sigma})$,  which preserves adjacency. Thus $\varphi$ induces an explicit isomorphism from $\Sigma$ to ${\hat \Sigma}$. By Lemma~\ref{lem:mixed-dih2}~(1), $\varphi$ also induces a permutational isomorphism between the $G$-actions on $\Sigma$ and ${\hat \Sigma}$. Thus it suffices to prove that part (a) holds for ${\hat \Sigma}$.

Recall that $N=H'$. By Lemma~\ref{lem:mixed-dih2}~(2), ${\hat \Sigma}$ is an $N$-normal cover of ${\hat \Sigma}_N\cong \K_{2^n,2^n}$.
Also, since $C(H,X,Y)$ is $G$-edge-transitive, it follows from Lemma~\ref{lem:mixed-dih2}~(4) that ${\hat \Sigma}$ is $(G,2)$-arc-transitive. It remains to prove that ${\hat \Sigma}_N\cong\K_{2^n,2^n}$ is $G/N$-edge-affine. By
Lemma~\ref{lem:mixed-dih2}~(1), $H$ acts regularly on the edge set of ${\hat \Sigma}$, and hence by Lemma~\ref{quot}~(4), $H/N$ acts regularly on the edge set of ${\hat \Sigma}_N\cong\K_{2^n,2^n}$. Since $H$ is an $n$-dimensional mixed dihedral group, we have $H/N\cong C_2^{2n}$ (Definition~\ref{mixed-dihedral-group}) and $H/N\unlhd G/N$  (since $H\unlhd G$), and hence ${\hat \Sigma}_N$ is $G$-edge-affine. Thus part (b) implies part (a).

Now assume that part (a) of Theorem~\ref{iff-characterization} holds, that is,  $\Sigma$ is a $(G,2)$-arc-transitive $N$-normal cover of $\K_{2^n,2^n}$, and  $\K_{2^n,2^n}$ is $G/N$-edge-affine. Thus in particular $\Sigma$ and $\K_{2^n,2^n}$ have the same valency, namely $2^n$. By Lemma~\ref{quot}~(2), $N$ is the kernel of the $G$-action on $\Sigma_N\cong K_{2^n,2^n}$, and hence $G/N\leq \Aut(\Sigma_N)$.
Since $\K_{2^n,2^n}$ is $G/N$-edge-affine, there exists $N<H\unlhd G$ such that $H/N\cong C_{2}^{2n}$ and $H/N$ is regular on the edge set of $\Sigma_N$. By Lemma~\ref{quot}~(4), $H$ is also regular on the edge set of $\Sigma$. Since $H'$ is characteristic in $H$ and $H\unlhd G$, we have $H'\unlhd G$. Consider the $H'$-normal quotient $\Sigma_{H'}$. Note that $\Sigma$ is bipartite, say with bipartition $V(\Sigma)= O_1\cup O_2$, and $O_1, O_2$ are the two $H$-orbits in $V(\Sigma)$.
Since $H/N\cong C_{2}^{2n}$ the derived subgroup $H'\leq N$, and since $N$ is intransitive on both $O_1$ and $O_2$, it follows that also $H'$ is intransitive on both $O_1$ and $O_2$.
Then since $G$ is $2$-arc-transitive on $\Sigma$, it follows from Lemma~\ref{quot} parts~(1) and~(4) that $\Sigma$ is an $H'$-normal cover of $\Sigma_{H'}$ (so $\Sigma_{H'}$ has valency $2^n$), and $H/H'$ is regular on the edge set of $\Sigma_{H'}$.

The group $H/H'$ is abelian and edge transitive on $\Sigma_{H'}$, and so by Lemma~\ref{lem:abelian-edge-regular},  $\Sigma_{H'}\cong \K_{2^n,2^n}$ (since $\Sigma_{H'}$ has valency $2^n\geq 4$). In particular, $|H/H'|= |E(\Sigma_{H'})|={2}^{2n}$. Since $H'\leq N<H$, it follows that $H'=N$ and $H/H'=H/N\cong C_{2}^{2n}$.
Since $H$ is regular on $E(\Sigma)$, the line graph $\mathcal{L}(\Sigma)$ of $\Sigma$ is a Cayley graph of $H$, say $\mathcal{L}(\Sigma)=\Cay(H,S)$. Let $\{u,v\}$ be an edge of $\Sigma$. Then
\[
S=\{h\in H:\ \{u,v\}^h\ {\rm is\ incident\ with}\ \{u,v\}\}.
\]
Now the set of edges of $\Sigma$ incident with $\{u,v\}$ is \[
F=\{\{u,x\}, \{v,y\}\mid v\neq x\in\Sigma(u), u\neq y\in\Sigma(v)\}.
\]
Let $X=H_u$ and $Y=H_v$. Then $XH'$ is the subgroup of $H$ stabilising the $H'$-orbit in $V(\Sigma)$ containing $u$. Since $H/H'\cong C_2^{2n}$ is regular on the set of edges of $\Sigma_{H'}\cong\K_{2^n,2^n}$, we have $XH'/H'\cong C_2^n$. Similarly,  $YH'/H'\cong C_2^n$. By Lemma~\ref{quot}~(2), $H'$ is semiregular on $V(\Sigma)$, so $X\cong XH'/H'\cong C_2^n$ and $Y\cong YH'/H'\cong C_2^n$. Again, since $H$ is regular on the edge set of $\Sigma$, it follows that $X$ and $Y$ act regularly on $\Sigma(u)$ and $\Sigma(v)$, respectively. Thus, $S=(X\cup Y)\setminus\{1\}$. Since $\Sigma$ is connected, $\mathcal{L}(\Sigma)$ is also connected, and so $H=\lg S\rg=\lg X, Y\rg$. It follows from Definition~\ref{mixed-dihedral-group} that $H$ is an $n$-dimensional mixed dihedral group relative to $X$ and $Y$, and from Definition~\ref{mixed-dihedrant} that $\mathcal{L}(\Sigma)=C(H,X,Y)$. Since $G$ is $2$-arc-transitive on $\Sigma$, it follows that $G_{\{u,v\}}$ is transitive on $F$, and so $\mathcal{L}(\Sigma)$ is $G$-edge-transitive. Thus part (b) holds. This completes the proof of Theorem~\ref{iff-characterization}.
\hfill\qed

\section{Construction of $2$-arc-transitive covers of $\K_{2^n,2^n}$}\label{constructionSection}

In this section, we apply Theorem~\ref{iff-characterization} to construct $2$-arc-transitive normal covers, of $2$-power order, of $\K_{2^n,2^n}$ and prove Theorem~\ref{2-dist-tran}. First we define in Definition~\ref{concreteConstruction} the group $\I(n)$ we will use in the construction. While it is  not obvious, the group $\I(n)$ is isomorphic to the group $\II(n)$ in Definition~\ref{constructionI}, a fact we prove in Lemma~\ref{CCpresentation}. The more explicit definition of multiplication for $\I(n)$ will help in our graph construction and analyses.


\begin{defi}\label{concreteConstruction}
Let $n\geq 2$, let $X$ and $Y$ be $n$-dimensional vector spaces over $\F_2$, and let $x_1,x_2\in X$, $y_1,y_2\in Y$ and $A_1,A_2\in X\otimes Y$. Define $\I(n)$ to be the
set $X\oplus Y\oplus (X\otimes Y)$ with multiplication defined as follows: for $g_1=x_1+y_1+A_1$ and $g_2=x_2+y_2+A_2$ in $\I(n)$, given by
\begin{equation}\label{eqnMultiplication}
    g_1g_2=g_1+g_2+x_2\otimes y_1,
\end{equation}
where each addition occurs in $X\oplus Y\oplus (X\otimes Y)$, considered as a vector space over $\F_2$.
\end{defi}

More explicitly, $g_1g_2$ is the element $(x_1+x_2)+(y_1+y_2) +(A_1+A_2+x_2\otimes y_1)$ of $\I(n)$. In particular we denote by $0$ the element $x+y+A\in \I(n)$ with each of $x, y, A$ equal to the zero vector of the corresponding space. It turns out that $\I(n)$ with this multiplication is a group.

\begin{lem}\label{CCgivesgroup}
 Let $X$, $Y$ and $\I(n)$ be as in Definition~\ref{concreteConstruction}. Then $\I(n)$ is a group of order $2^{n^2+2n}$ with identity $0$. Furthermore, for $g=x+y+A\in \I(n)$, with $x\in X,y\in Y,A\in X\otimes Y$, the inverse is
 \begin{equation}\label{eqnInverse}
     g^{-1}=g+x\otimes y.
 \end{equation}
\end{lem}

\f\demo
By definition $|\I(n)|=2^{n^2+2n}$.
 First, we show that the multiplication is associative. For $i=1,2,3$, let $g_i=x_i+y_i+A_i\in\I(n)$, where $x_i\in X,y_i\in Y,A_i\in X\otimes Y$. Then, applying the multiplication defined in Equation~(\ref{eqnMultiplication}),
 \begin{align*}
  (g_1g_2)g_3&=((x_1+x_2)+(y_1+y_2) +(A_1+A_2+x_2\otimes y_1))g_3 \\
 &= (x_1+x_2+x_3)+(y_1+y_2+y_3) +(A_1+A_2+A_3+x_2\otimes y_1+x_3\otimes (y_1+y_2))\\
  &=g_1((x_2+x_3)+(y_2+y_3) +(A_2+A_3+x_3\otimes y_2)) \\
  &=g_1(g_2g_3).
 \end{align*}
 A direct computation using \eqref{eqnMultiplication} yields  $g0=g=0g$, so $0$ is the identity of $\I(n)$. Finally, for $g=x+y+A\in\I(n)$,  using \eqref{eqnMultiplication} we have
 \[
  g(g+x\otimes y)=g + (g+x\otimes y) + x\otimes y = 2g+2(x\otimes y)=0
 \]
 and similarly $(g+x\otimes y)g=0$.
 Hence $g^{-1}=g+x\otimes y$. In particular $\I(n)$ is a group.
\hfill\qed

Next we derive some relevant properties of $\I(n)$.

\begin{lem}\label{CCismixedDihedral}
 Let $X$, $Y$ and $\I(n)$ be as in Definition~\ref{concreteConstruction}. Then the following hold.
 \begin{enumerate}
     \item [{\rm (1)}] As subgroups of $\I(n)$, we have $X\cong Y\cong C_2^n$ and $X\otimes Y\cong C_2^{n^2}$.
     \item [{\rm (2)}] $\I(n)=\langle X,Y\rangle$.
     \item [{\rm (3)}] The derived subgroup $\I(n)'=X\otimes Y$, and $\I(n)'$ is equal to the centre $Z(\I(n))$.
     \item [{\rm (4)}] $\I(n)/\I(n)'\cong C_2^{2n}$ and $|\I(n)|=2^{n^2+2n}$.
 \end{enumerate}
 In particular, $\I(n)$ is an $n$-dimensional mixed dihedral group with respect to $X$ and $Y$.
\end{lem}

\f\demo
 Throughout the proof, for each $i=1,2$, we let $x_i\in X$, $y_i\in Y$, $A_i\in X\otimes Y$ and $g_i\in\I(n)$ such that $g_i=x_i+y_i+A_i$.

 (1) By Equation~(\ref{eqnMultiplication}), $x_1x_2=x_1+x_2=x_2x_1$ and $x_1^2=0$. Since $|X|=2^n$, we have $X\cong C_2^n$. By a similar argument, $Y\cong C_2^n$. Also, $A_1A_2=A_1+A_2=A_2A_1$, and the product $A_1A_1=2A_1=0$ and $|X\otimes Y|=2^{n^2}$. Hence $X\otimes Y\cong C_2^{n^2}$ and part (1) is proved.

(2) A typical element of $\I(n)$ has the form $g_1=x_1+y_1+A_1$. By Equation~(\ref{eqnMultiplication}), the product $x_1y_1A_1=x_1+y_1+A_1=g_1$ and so, in order to prove that $\I(n)=\langle X,Y\rangle$ it suffices to prove that $A_1$ lies in the subgroup $\langle X,Y\rangle$. As $\F_2$-vector spaces, let $\{e_1,\ldots,e_n\}$ be a basis for $X$ and $\{f_1,\ldots,f_n\}$ be a basis for $Y$ so that $\{e_i\otimes f_j\mid 1\leq i,j\leq n\}$ is a basis for $X\otimes Y$. Since $A_1\in X\otimes Y$, there exist $c_{ij}\in\{0,1\}$, for $1\leq i,j\leq n$, such that
 \[
  A_1=\sum_{i=1}^n\sum_{j=1}^n c_{ij} (e_i\otimes f_j).
 \]
 By part (1), $X\otimes Y$ is an abelian subgroup of $\I(n)$, and hence, by the definition of multiplication \eqref{eqnMultiplication}, we can rewrite the above summation as
 \[
  A_1=\prod_{i=1}^n\prod_{j=1}^n (e_i\otimes f_j)^{c_{ij}}.
 \]
 Again using \eqref{eqnMultiplication},  $(e_if_j)^2=(e_i+f_j)^2=e_i\otimes f_j$, and hence $e_i\otimes f_j\in\langle X, Y\rangle$, for each $i,j$. It then follows from the product expression for $A_1$ displayed above that $A_1\in \langle X, Y\rangle$. Therefore $\I(n)=\langle X, Y\rangle$, and  part (2) is proved.

 (3)  Recalling the form for the inverse of an  element of $\I(n)$ in Equation~(\ref{eqnInverse}), we have
 \begin{align*}
  [g_1,g_2]&=g_1^{-1}g_2^{-1}g_1g_2=(g_1+x_1\otimes y_1)(g_2+x_2\otimes y_2)g_1g_2\\
  &=(g_1+g_2+x_1\otimes y_1+x_2\otimes y_2+x_2\otimes y_1)(g_1+g_2+x_2\otimes y_1)\\
  &=x_1\otimes y_1+x_2\otimes y_2+(x_1+x_2)\otimes (y_1+y_2)\\
  &=x_1\otimes y_2+x_2\otimes y_1.
 \end{align*}
 Thus $\I(n)'\leq X\otimes Y$. Moreover, in terms of the basis for $X\otimes Y$ introduced in (2) above, we have $[e_i,f_j]=(e_i+f_j)^2=e_i\otimes f_j$ for each $i,j$. Hence $X\otimes Y$ is generated by the commutators $[e_i,f_j]$, for $1\leq i,j\leq n$, so $\I(n)'= X\otimes Y$. In the displayed equation above, choosing $g_2=A_2\in X\otimes Y=\I(n)'$ leads to $[g_1,g_2]=0$ for all $g_1\in\I(n)$, and hence $\I(n)'\leq Z(\I(n))$.  On the other hand, if $g_2\in \I(n)\setminus \I(n)'$, then at least one of $x_2, y_2\ne 0$. Without loss of generality, $x_2=\sum_{i=1}^n a_ie_i$ with, say, $a_\ell=1$. Then
 choosing $g_1=y_1=f_\ell$, we have $[g_1,g_2]=x_2\times y_1 = \sum_{i=1}^n a_i e_i\otimes f_\ell$, which is nonzero since $a_\ell\ne0$. Thus $g_2\not\in Z(\I(n))$, and hence $\I(n)'=Z(\I(n))$, proving part (3).

(4) Finally, working in the quotient $\I(n)/\I(n)'$, $g_1 (X\otimes Y)g_2 (X\otimes Y)=(x_1+x_2+y_1+y_2)(X\otimes Y)=g_2 (X\otimes Y)g_1 (X\otimes Y)$, and $(g_1 (X\otimes Y))^2=X\otimes Y$, so $\I(n)/\I(n)'$ is an elementary abelian $2$-group. By Lemma~\ref{CCgivesgroup}, $|\I(n)|=2^{n^2+2n}$, and by part (3), $|\I(n)'|=2^{n^2}$, so part (4) follows. The last assertion follows from Definition~\ref{mixed-dihedral-group}, completing the proof.
\hfill\qed

Now we prove the isomorphism of the groups $\I(n)$ and $\II(n)$, and  we note that the assertions of Theorem~\ref{2-dist-tran}~(1) then follow from Lemmas~\ref{CCismixedDihedral} and~\ref{CCpresentation}.

\begin{lem}\label{CCpresentation}
 For each $n\geq2$, the group $\I(n)$  in Definition~\ref{concreteConstruction} is isomorphic to the group $\II(n)$ given by the presentation in Definition~\ref{constructionI}.
\end{lem}

\f\demo
 We claim that $\I(n)$, as defined in Definition~\ref{concreteConstruction}, satisfies all the relations in Definition~\ref{constructionI}. Let $e_1,\ldots,e_n$ be a basis for $X$ and $f_1,\ldots,f_n$ be a basis for $Y$. By Lemma~\ref{CCismixedDihedral} (2) we have that $\I(n)=\langle e_1,\ldots,e_n,f_1,\ldots,f_n\rangle$. Since the multiplication defined in Equation~(\ref{eqnMultiplication}) reduces to vector addition when restricted to elements of $X$, or, respectively, when restricted to elements of $Y$, we have $e_i^2=0$, $f_i^2=0$, $[e_i,e_j]=0$ and $[f_i,f_j]=0$ for all $i,j\in\{1,\ldots,n\}$. Also, for each $i,j\in\{1,\ldots,n\}$ we have $[e_i,f_j]=(e_i+f_j)^2=e_i\otimes f_j$. It then follows from Equation~(\ref{eqnMultiplication}) that $[e_i,f_j]^2=0$ and that $[e_i,f_j]$ commutes with $e_k$ and $f_k$ for all $k\in\{1,\ldots,n\}$. This proves the claim.

 It follows that $\I(n)$ is isomorphic to a quotient of the group $\II(n)$  in Definition~\ref{constructionI}. %
 Moreover, it follows from the relations in Definition~\ref{constructionI} that every element of the group $\II(n)$ may be written as a product
 \[
  \left(\prod_{i=1}^n x_i^{a_i}\right) \left(\prod_{i=1}^n y_i^{b_i}\right) \left(\prod_{i=1}^n\prod_{j=1}^n [x_i,y_j]^{c_{ij}}\right),
 \]
 for some $a_i,b_i,c_{ij}\in\{0,1\}$, where $i,j\in\{1,\ldots,n\}$. The number of distinct tuples of parameters $a_i, b_i, c_{ij}$ is $2^{n^2+2n}$, and hence $\II(n)$ has order at most $2^{n^2+2n}$. Since by Lemma~\ref{CCismixedDihedral}~(4), $|\I(n)|=2^{n^2+2n}$, we conclude that $|\I(n)|=|\II(n)|$ and hence that $\I(n)\cong\II(n)$, completing the proof.
 \hfill\qed

Our next task is to determine the subgroup of $\Aut(\I(n))$ that leaves $X\cup Y$ invariant.

\begin{lem}\label{CCautomorphsimGroup}
 Let $X$, $Y$, $\I(n)$ be as in Definition~\ref{concreteConstruction} and let $S=(X\cup Y)\setminus\{0\}$. Then the subgroup $\Aut(\I(n),S)$ of $\Aut(\I(n))$ stabilising $S$ is $(\Aut(X)\times\Aut(Y)):C_2\cong \GL_n(2)\wr S_2$.
\end{lem}

\f\demo
 First, by Lemma~\ref{CCismixedDihedral}, $X\cong Y\cong C_2^n$, and hence $\Aut(X)\cong\Aut(Y)\cong\GL_n(2)$. As $\F_2$-vector spaces, let $\{e_1,\ldots,e_n\}$ be a basis for $X$ and let $\{f_1,\ldots,f_n\}$ be a basis for $Y$ so that $\{e_i\otimes f_j\mid 1\leq i,j\leq n\}$ is a basis for $X\otimes Y$. Let
 \[
  g_1=x_1+y_1+\sum_{i=1}^n\sum_{j=1}^n a_{ij} (e_i\otimes f_j) \quad\text{and}\quad g_2=x_2+y_2+\sum_{i=1}^n\sum_{j=1}^n b_{ij} (e_i\otimes f_j),
 \]
 where $x_1,x_2\in X$, $y_1,y_2\in Y$ and $a_{ij},b_{ij}\in\F_2$ for $1\leq i,j\leq n$. For $\phi\in\Aut(X)$ we define
 \[
  g_1^\phi=x_1^\phi+y_1+\sum_{i=1}^n\sum_{j=1}^n a_{ij} (e_i^\phi\otimes f_j).
 \]
 Since $\phi$ is an invertible linear transformation of $X$, it follows that $\phi$ defines a bijection on $X\oplus Y\oplus (X\otimes Y)$. Also,
 \begin{align*}
  g_1^\phi g_2^\phi&=\left(x_1^\phi+y_1+\sum_{i=1}^n\sum_{j=1}^n a_{ij} (e_i^\phi\otimes f_j)\right)\left(x_2^\phi+y_2+\sum_{i=1}^n\sum_{j=1}^n b_{ij} (e_i^\phi\otimes f_j)\right)\\
  &=x_1^\phi+x_2^\phi+y_1+y_2+\sum_{i=1}^n\sum_{j=1}^n (a_{ij}+b_{ij}) (e_i^\phi\otimes f_j)+x_2^\phi\otimes y_1,
 \end{align*}
 while
 \begin{align*}
  (g_1 g_2)^\phi&=\left(x_1+x_2+y_1+y_2+\sum_{i=1}^n\sum_{j=1}^n (a_{ij}+b_{ij}) (e_i\otimes f_j)+x_2\otimes y_1\right)^\phi\\
  &=x_1^\phi+x_2^\phi+y_1+y_2+\sum_{i=1}^n\sum_{j=1}^n (a_{ij}+b_{ij}) (e_i^\phi\otimes f_j)+x_2^\phi\otimes y_1.
 \end{align*}
 Thus $g_1^\phi g_2^\phi=(g_1 g_2)^\phi$ and so $\phi\in\Aut(\I(n))$ and $\phi$ leaves $X\cup Y$ invariant. Similarly if $\phi\in\Aut(Y)$ then the map
 \[
  g_1^\phi=x_1+y_1^\phi+\sum_{i=1}^n\sum_{j=1}^n a_{ij} (e_i\otimes f_j^\phi)
 \]
 defines an element of $\Aut(\I(n))$ and $\phi$ leaving $X\cup Y$ invariant. Thus $\Aut(X)\times\Aut(Y)\leq \Aut(\I(n),S)$.
 Further, for
 \[
  g=\sum_{i=1}^n a_ie_i + \sum_{i=1}^n b_if_i + \sum_{i=1}^n\sum_{j=1}^n c_{ij} (e_i\otimes f_j)\  \in\I(n),
 \]
 define a map $\delta:\I(n)\to \I(n)$ by
 \[
  g^\delta := \sum_{i=1}^n b_ie_i + \sum_{i=1}^n a_if_i + \sum_{i=1}^n\sum_{j=1}^n (a_jb_i+c_{ji}) (e_i\otimes f_j).
 \]
 Observe that $\delta$ is the composition of the map that interchanges $g$ and $g^{-1}$ (see Equation~(\ref{eqnInverse})) with the map that interchanges $e_i$ and $f_i$, for $1\leq i\leq n$. It follows that $\delta$ is a bijection, and also that $S^\delta=S$. To see that $\delta$ is an automorphism consider $g$ as above and
 \[
  g'=\sum_{i=1}^n a_i'e_i + \sum_{i=1}^n b_i'f_i + \sum_{i=1}^n\sum_{j=1}^n c_{ij}' (e_i\otimes f_j),
 \]
 and compute as follows:
 \begin{align*}
  g^\delta (g')^\delta&=\left( \sum_{i=1}^n b_ie_i + \sum_{i=1}^n a_if_i + \sum_{i=1}^n\sum_{j=1}^n (a_jb_i+c_{ji}) (e_i\otimes f_j)\right)\\
  &\phantom{white}\left( \sum_{i=1}^n b_i'e_i + \sum_{i=1}^n a_i'f_i + \sum_{i=1}^n\sum_{j=1}^n (a_j'b_i'+c_{ji}') (e_i\otimes f_j)\right)\\
  &= \sum_{i=1}^n (b_i+b_i')e_i + \sum_{i=1}^n (a_i+a_i')f_i + \sum_{i=1}^n\sum_{j=1}^n (a_jb_i+a_jb_i'+a_j'b_i'+c_{ji}+c_{ji}') (e_i\otimes f_j).
 \end{align*}
 and
 \begin{align*}
  (gg')^\delta&= \left(\sum_{i=1}^n (a_i+a_i')e_i + \sum_{i=1}^n (b_i+b_i')f_i + \sum_{i=1}^n\sum_{j=1}^n (a_i'b_j+c_{ij}+c_{ij}') (e_i\otimes f_j)\right)^\delta\\
  &=\sum_{i=1}^n (b_i+b_i')e_i + \sum_{i=1}^n (a_i+a_i')f_i + \sum_{i=1}^n\sum_{j=1}^n (a_jb_i+a_jb_i'+a_j'b_i'+c_{ji}+c_{ji}') (e_i\otimes f_j).
 \end{align*}
 Thus $\delta\in \Aut(\I(n),S)$, and hence we have proved that $\Aut(\I(n),S)$ contains $(\Aut(X)\times\Aut(Y)):C_2$. However, by Lemma~\ref{lem:mixed-dih2}~(3),  $\Aut(\I(n),S)$ is a subgroup of $(\Aut(X)\times\Aut(Y)):C_2$, and hence equality holds, completing the proof.
\hfill\qed

Now we investigate the graphs $\G(n):=C(\I(n),X,Y)$ and $\Sigma(n):=\Sigma(\I(n),X,Y)$ for $\I(n)$,  as defined in Definition~\ref{mixed-dihedrant}. Figure~\ref{distdiaggamm22} shows a distance diagram for the smallest of these graphs $\G(2)$. We note that Theorem~\ref{2-dist-tran}~(2) follows from Lemma~\ref{fullAutomorphism}~(2).

\begin{lem}\label{fullAutomorphism}
Let $X$, $Y$, $\I(n)$ be as in Definition~\ref{concreteConstruction}, let $\G(n)=C(\I(n),X,Y)$, and $\Sigma(n)=\Sigma(\I(n), X, Y)$, as in Definition~{\rm \ref{mixed-dihedrant}}, and let $S=(X\cup Y)\setminus\{0\}$ and $N=\I(n)'$. Then
\begin{enumerate}
    \item[{\rm (1)}] $\Aut(\G(n))=\Aut(\Sigma(n))=\I(n): \Aut(\I(n), S)\cong \I(n): (\GL_n(2)\wr S_2)$.

    \item[{\rm (2)}] $\Sigma(n)$ is a $2$-arc-transitive graph of order $2^{n^2+n+1}$, and an $N$-normal cover of $\K_{2^n,2^n}$, which is $\Aut(\Sigma(n))/N$-edge-affine, and in particular is $\I(n)/N$-edge-affine.
    \end{enumerate}
\end{lem}



\begin{rem}\label{rem-magma}
{\rm In the following proof we make use of a Magma~\cite{BCP} computation in the cases $n=2$ or $3$, which we now describe. First, the group $\II(n)$, with $n=2$ or $3$, is input in the category GrpFP via the presentation given in Definition~\ref{constructionI}. Next, the {\tt pQuotient} command is used to construct the largest $2$-quotient $H$ of $\II(n)$ having lower exponent-$2$ class at most 100 as group in the category GrpPC. Comparing the orders of these groups, we find $|\II(n)|=|H|$, so that $\II(n)\cong H$. We then construct a graph isomorphic to $\G(n)$ as a Cayley graph on $H$. Computing the order of the full automorphism of this graph then shows that $\Aut(\G(n))=A=\I(n): \Aut(\I(n), S)$. We have made available the Magma programs in the Appendix of this paper.}
\end{rem}


\f\demo Write $\G=\G(n)$ and $\Sigma=\Sigma(n)$. By Lemma~\ref{CCautomorphsimGroup}, the automorphism group $\Aut(\G)$ contains as a subgroup
\[
A:= \I(n): \Aut(\I(n), S)\cong \I(n): (\GL(n,2)\wr S_2).
\]
Also $\G$ is $A$-edge-transitive, and $\I(n)$ acts regularly on $E(\Sigma)$ and has two orbits on $V(\Sigma)$, by Lemma~\ref{lem:mixed-dih2}~(1).
By Lemma~\ref{lem:prop-mixed-dih}~(4), the map $\varphi:z\to \{Xz,Yz\}$ defines an isomorphism from $\G$ to the line graph of $\Sigma$, and  by Lemma~\ref{lem:prop-mixed-dih}~(5), $\Aut(
 \G)=\Aut(\Sigma)$. By Lemma~\ref{CCismixedDihedral}~(3), the $2$-group $N:=\I(n)'=Z(\I(n))$ is both the derived subgroup and the centre of $\I(n)$.
Furthermore, by Lemma~\ref{lem:mixed-dih2}~(2), the $N$-normal quotient $\Sigma_N$ is isomorphic to $\K_{2^n,2^n}$. Since (as noted above) $A$ is edge-transitive on $\G(n)$, it follows from Lemma~\ref{lem:mixed-dih2}~(4) that $A$ is $2$-arc-transitive on $\Sigma$, and then from Lemma~\ref{quot}~(3) that $A/N$ is $2$-arc-transitive on  $\Sigma_N$.
By Lemma~\ref{lem:mixed-dih2}~(2) $\Sigma$ is an $N$-normal cover of $\K_{2^n,2^n}$. Hence $\Sigma$ has order $2^{n+1}\cdot |N|$. By Lemma~\ref{CCismixedDihedral}~(4), we have $|N|=2^{n^2}$. Thus $\Sigma$ has order $2^{n^2+n+1}$.
We now apply Lemma~\ref{sufficient} to $\Sigma$ with $2$-arc-transitive group $A$. It follows from Lemma~\ref{CCautomorphsimGroup}  that the kernel of the action of the stabiliser $A_X$ on $\{Yx : x\in X\}$ contains $\Aut(Y)$ and in particular is nontrivial, and hence, the last assertion of Lemma~\ref{sufficient} implies that $N\unlhd \Aut(\Sigma)$.

For any $T$ such that $N\leq T\leq  \Aut(\Sigma)$,
let $T^+/N$ denote the subgroup of $T/N$ stabilising both parts of the bipartition of $\Sigma_{N} \cong\K_{2^n,2^n}$.
By Lemma~\ref{quot}~(2), the kernel of the action of
$\Aut(\Sigma)$ on $\Sigma_N$ is $N$ so $\I(n)/N\leq \Aut(\Sigma)^+/N\leq \Aut(\Sigma_N)$, and also,  by  Lemma~\ref{CCismixedDihedral}~(4), $\I(n)/N\cong C_2^{2n}$.
By Lemma~\ref{lem:mixed-dih2}~(1), $\I(n)$ is regular on $E(\Sigma)$, and hence $\I(n)/N$ is transitive on $E(\Sigma_{N})$; moreover, since $\I(n)/N$ is abelian the latter action is regular.
Since $\I(n)$ has two orbits on $V(\Sigma)$, it follows that $\I(n)/N$ has two orbits on $V(\Sigma_N)$, and since $\I(n)/N$ is edge-regular, it follows that $\Sigma_N$ is $A/N$-edge-affine. Note that $\Aut(\Sigma_N)=S_{2^n}\wr S_2$. If $\Aut(\Sigma)^+/N$ has socle  $C_2^{2n}$, then
$A/N\leq \Aut(\Sigma)/N\leq C_2^{2n} : (\GL(n,2)\wr S_2) \cong A/N$, and it follows that  $\Aut(\Sigma)/N = A/N$ and $\Aut(\Sigma)=A$. In this case $\Sigma_N$ is $\Aut(\Sigma)/N$-edge-affine and both parts of the lemma are proved in this case.

For the small values $n=2$ or $3$, we check, using Magma~\cite{BCP}, that $A=\Aut(\G(n))$, and hence the lemma is proved for  $n=2$ or $3$ (see Remark~\ref{rem-magma} for a description of these computations).  Now assume that $n\geq 4$, and assume, for a contradiction, that the socle of $\Aut(\Sigma_N)^+$ is not $C_2^{2n}$.
We will apply Proposition~\ref{basic-auto} to $\Sigma_{N} \cong\K_{2^n,2^n}$ and $\Aut(\Sigma)/N$. Since $n\geq4$, case  (4) of Proposition~\ref{basic-auto} does not hold.
In case~(2) of Proposition~\ref{basic-auto}, the group $\Aut(\Sigma)^+/N$ does not contain a subgroup isomorphic to $C_2^{2n}$, so this case does not hold either. Thus case (1) of Proposition~\ref{basic-auto} holds, so we have normal subgroups $M, T_1, T_2$ of $\Aut(\Sigma)^+$, all containing $N$, such that $M/N=\soc(\Aut(\Sigma)/N)=T_1/N\times T_2/N\cong A_{2^n}\times A_{2^n}$.
We showed above that $\I(n)/N$ is regular on $E(\Sigma_{N})$ and has two vertex-orbits.
It follows that  $\I(n)/N=(\I(n)/N)^{U}\times (\I(n)/N)^{W}$, where $U,W$ are the two parts of the bipartition of $\Sigma_{N}$, and $(\I(n)/N)^U$ and $(\I(n)/N)^{W}$ are the permutation groups induced by $\I(n)/N$ on $U$ and $W$, respectively.
Moreover, $(\I(n)/N)^U, (\I(n)/N)^W$ is regular on $U, W$, respectively. Since $n\geq4$, the permutation induced on $U$ by a nontrivial element of $(\I(n)/N)^U$ is a product of $2^{n-1}$ cycles of length $2$, and hence is an even permutation of $U$. Similarly every element of $(\I(n)/N)^W$ induces an even permutation of $W$. Thus $\I(n)/N\leq M/N\cong A_{2^n}\times A_{2^n}$, and the normaliser of $\I(n)/N$ in $M/N$ is isomorphic to $C_2^{2n}: (\GL(n,2)\times \GL(n,2))$. It follows that the normaliser of $\I(n)/N$ in $M/N$ is $A^+/N=(\I(n):(\Aut(X)\times\Aut(Y)))/N$.
We may assume that $\Aut(X)N/N\leq T_1/N$ and $\Aut(Y)N/N\leq T_2/N$. Now  $N\cong C_2^{n^2}$,  and so $N\leq C_M(N)$ and $M/C_M(N)\leq\Aut(N)\cong\GL(n^2,2)$. We saw in the proof of Lemma~\ref{CCautomorphsimGroup} that
both $\Aut(X)$ and $\Aut(Y)$ act nontrivially and faithfully on $N$, so in particular neither $T_1/N$ nor $T_2/N$ is contained in $C_M(N)$. Hence $C_M(N)/N$ is a proper normal subgroup of $M/N\cong A_{2^n}\times A_{2^n}$, containing neither simple direct factor. It follows that $C_M(N)/N=1$. Hence $M/C_M(N)=M/N\cong A_{2^n}\times A_{2^n}$. Since $M/C_M(N)\leq\GL(n^2,2)$ and since $2^n\geq 16$, it follows from \cite[Proposition~5.3.2]{Kleidman-Liebeck} that $n^2\geq 2^n-2$, which  implies that $n=4$. However $|M/N|=|A_{16}|^2$ is divisible by $13^2$, while $13^2$ does not divide $|\GL(16,2)|$ (as the least $j$ such that $13$ divides $2^j-1$ is $j=12$), and hence $13^2$ does not divide $|M/C_M(N)|$. This contradiction completes the proof.
\hfill\qed

We next prove various symmetry properties of $\G(n)$ which complete the proof of Theorem~\ref{2-dist-tran}~(3).

\begin{lem}\label{answer-prob-3-new}
Let $X$, $Y$, $\I(n)$ be as in Definition~\ref{concreteConstruction}, let $\G(n)=C(\I(n),X,Y)$, as in Definition~{\rm \ref{mixed-dihedrant}}, and let $S=(X\cup Y)\setminus\{0\}$. Then $\G(n)$ is a $2$-geodesic-transitive normal Cayley graph; moreover $\G(n)$ is $2$-distance-transitive, but is neither $3$-distance-transitive (and in particular not distance-transitive), nor $2$-arc-transitive.
\end{lem}

\f\demo By Definition~\ref{mixed-dihedrant}, $\G:=\G(n)$ is a Cayley graph, and by Lemma~\ref{fullAutomorphism}, $\Aut(\G)=\I(n):\Aut(\I(n),S)\cong \I(n):(\GL_n(2)\wr S_2)$. Thus $\G$ is a normal Cayley graph. Let $G=\Aut(\G)$. Since $\I(n)$ acts transitively on $V(\G)$, to prove the result it is sufficient to consider, for the vertex $0\in\I(n)=V(\G)$, the action of $G_0=\Aut(\I(n),S)$ on the vertices at distance one, two and three from $0$. The set of vertices that are distance one from $0$ is precisely $\G(0)=S$, on which $G_0$ acts transitively, so $\G$ is 1-arc-transitive.

We now consider the set $\Gamma_2(0)$ of vertices at distance two from $0$. Let $X_0=X\setminus\{0\}$, let $Y_0=Y\setminus\{0\}$, let $x,x'$ be distinct elements of $X_0$ and let $y,y'$ be distinct elements of $Y_0$. Then, recalling that multiplication is as in Equation~(\ref{eqnMultiplication}), $xx'\in X$, $yy'\in Y$, $xy=x+y\in \I(n)\setminus (X\cup Y)$ and $yx=x+y+x\otimes y\in \I(n)\setminus (X\cup Y)$. Hence, the set of vertices at distance two from $0$ is
\[
 \Gamma_2(0)=\{x+y\mid x\in X_0,y\in Y_0\}\cup \{x+y+x\otimes y\mid x\in X_0,y\in Y_0\}.
\]
Let $e_1,\ldots,e_n$ be a basis for $X$, let $f_1,\ldots,f_n$ be a basis for $Y$, and let $\delta$ be as in the proof of Lemma~\ref{CCautomorphsimGroup}, that is, $\delta$ is the composition of the map interchanging $g$ and $g^{-1}$ (for $g\in\I(n)$, see Equation~(\ref{eqnInverse})) with the map interchanging $e_i$ and $f_i$, for $i\in\{1,\ldots,n\}$. Then, there exist $\sigma$ and $\sigma'$ in $\Aut(X)\times\Aut(Y)$ such that $(x+y)^\sigma=e_1+f_1$ and $(x+y+x\otimes y)^{\sigma'}=e_1+f_1+e_1\otimes f_1$. Since $(e_1+f_1)^\delta=e_1+f_1+e_1\otimes f_1$, we conclude that $(x+y)^{\sigma\delta\sigma'}=x+y+x\otimes y$. Hence $\G_2(0)$ is a $G_0$-orbit, and thus $\G$ is $2$-distance-transitive. Considering the right multiplication action by elements of $X\cup Y$ on various edges $\{0,s\}$, for $s\in S$, we see that
\[
\G(x) = \{0\}\cup \{x' : x'\in X_0, x'\ne x\} \cup \{yx : y\in Y_0\}
\]
and
\[
\G(y) = \{0\}\cup \{y' : y'\in Y_0, y'\ne y\} \cup \{xy : x\in X_0\}
\]
and hence the only vertex at distance one from both $0$ and $xy=x+y$ is $y$. It follows that $\G$ is $2$-geodesic-transitive.

We now determine the set $\Gamma_3(0)$ of vertices that are at distance three from $0$, each such vertex is of the form $sz$ for some $s\in S, z\in\G_2(0)$ (note that we obtain the edge $\{z,sz\}$  by right-multiplying the edge $\{0,s\}$ by the element $z\in \G_2(0)$). We have the following cases, where $x,x'$ are distinct elements of $X_0$, $y,y'$ are distinct elements of $Y_0$; we set  $x''=xx'$ and $y''=yy'$, and note that $x''\not\in\{0,x,x'\}$ and $y''\not\in\{0,y,y'\}$:
\begin{align*}
 x(x+y)&=y\in S,\\
 x(x+y+x\otimes y)&=y+x\otimes y\in\I(n)\setminus(X\cup Y\cup\G_2(0)),\\
 y(x+y)&=x+x\otimes y\in\I(n)\setminus(X\cup Y\cup\G_2(0)),\\
 y(x+y+x\otimes y)&=x\in S, \\
 x'(x+y)&=x''+y\in \G_2(0),\\
 x'(x+y+x\otimes y)&=x''+y+x\otimes y\in\I(n)\setminus(X\cup Y\cup\G_2(0)),\\
 y'(x+y)&=x+y''+x\otimes y'\in\I(n)\setminus(X\cup Y\cup\G_2(0)) \quad\text{and}\\
 y'(x+y+x\otimes y)&=x+y''+x\otimes y''\in \G_2(0).
\end{align*}
 Hence, the set of vertices at distance three from $0$ is
\begin{align*}
 \Gamma_3(0)=\{y+x\otimes y\mid x\in X_0,y\in Y_0\}&\cup \{x+y+x\otimes y'\mid x\in X_0;y,y'\in Y_0,y\neq y'\}\\
 \cup\{x+x\otimes y\mid x\in X_0,y\in Y_0\}&\cup\{x+y+x'\otimes y\mid x,x'\in X_0,x\neq x';y\in Y_0\}.
\end{align*}
Now, each of the sets $\{y+x\otimes y\mid x\in X_0,y\in Y_0\}$ and $\{x+x\otimes y\mid x\in X_0,y\in Y_0\}$ are $(\Aut(X)\times\Aut(Y))$-orbits, and these sets are interchanged by $\delta$. Thus
\[
 \{y+x\otimes y\mid x\in X_0,y\in Y_0\}\cup\{x+x\otimes y\mid x\in X_0,y\in Y_0\}
\]
is a $G_0$-orbit. Hence $G_0$ does not act transitively on $\G_3(0)$ and thus $\G$ is not $3$-distance-transitive. To see that $\G$ is not $2$-arc-transitive, we observe that $(x,0,x')$ and $(x,0,y)$ are $2$-arcs and that $\{x,x'\}$ is an edge, while $\{x,y\}$ is not, and hence there is no element of $\Aut(\G)$ mapping $(x,0,x')$ to $(x,0,y)$. This completes the proof.
\hfill\qed

Finally we formalise the proof of Theorem~\ref{2-dist-tran}.

\bigskip\noindent
\textbf{Proof of Theorem~\ref{2-dist-tran}}\quad  As noted above, Theorem~\ref{2-dist-tran}~(1) follows from Lemmas~\ref{CCismixedDihedral} and~\ref{CCpresentation}, while Theorem~\ref{2-dist-tran}~(2) follows from Lemma~\ref{fullAutomorphism}~(2). Finally Theorem~\ref{2-dist-tran}~(3) follows from Lemma~\ref{answer-prob-3-new}.





\section*{Appendix: Magma programs for proving Lemma~\ref{fullAutomorphism} in case $n=2$ or $3$.}

\f{\bf Appendix~1.}\ A function for constructing Cayley graphs:\smallskip

{\tt Cay:=function(G,S);

V:={g:g in G};

E:={{g,s*g}:g in G,s in S};

return Graph<V|E>;

end function;}\smallskip

\f{\bf Appendix~2.}\ The case $n=2$:\medskip

{Input the group $\II(2)$:}\smallskip

{\tt G<x1,x2,y1,y2>:=Group< x1,x2,y1,y2 | x1\textasciicircum 2,x2\textasciicircum 2,y1\textasciicircum 2,y2\textasciicircum 2,

(x1,x2)=(y1,y2)=1,

(x1,y1)\textasciicircum 2=(x1,y2)\textasciicircum 2=(x2,y1)\textasciicircum 2=(x2,y2)\textasciicircum 2=1,

((x1,y1),x1)=((x1,y1),x2)=((x1,y1),y1)=((x1,y1),y2)=1,

((x1,y2),x1)=((x1,y2),x2)=((x1,y2),y1)=((x1,y2),y2)=1,

((x2,y1),x1)=((x2,y1),x2)=((x2,y1),y1)=((x2,y1),y2)=1,

((x2,y2),x1)=((x2,y2),x2)=((x2,y2),y1)=((x2,y2),y2)=1>;}
\smallskip

{Construct the largest 2-quotient group of $\II(2)$ having lower exponent-$2$ class at most 100 as group in the category GrpPC:}
\smallskip

{\tt G2,q:=pQuotient(G,2,100)};\smallskip

{Order of $G2$ (The result shows that $|G2|=\II(2)$, and so $G2\cong \II(2)$):}\smallskip

{\tt FactoredOrder(G2);}\smallskip

Construct the graph $\Gamma(2)=\Cay(G2, S)$:
\smallskip

{\tt x1:=x1@q;
x2:=x2@q;
y1:=y1@q;
y2:=y2@q;

S:={x1,x2,x1*x2,y1,y2,y1*y2};

Gamma2:=Cay(G2,S);}
\smallskip

{Automorphism Group of $\Gamma(2)$ (The result shows that $|\Aut(\Gamma(2))|=|G(2)||\GL_2(2)\wr S_2|$):}\smallskip

{\tt A:=AutomorphismGroup(Gamma2);

$\sharp$A eq ($\sharp$G2)*$\sharp$(GL(2,2))*$\sharp$(GL(2,2))*2;}

\medskip
\f{\bf Appendix~3.}\ The case $n=3$:\medskip

Input group $\II(3)$:\smallskip

{\tt G<x1,x2,x3,y1,y2,y3>:=Group<x1,x2,x3,y1,y2,y3 | x1\textasciicircum 2, x2\textasciicircum 2, x3\textasciicircum 2, y1\textasciicircum 2,

y2\textasciicircum 2, y3\textasciicircum 2, (x1,x2)=(x1,x3)=(x2,x3)=(y1,y2)=(y1,y3)=(y2,y3)=1,

(x1,y1)\textasciicircum 2=(x1,y2)\textasciicircum 2=(x1,y3)\textasciicircum 2=(x2,y1)\textasciicircum 2=(x2,y2)\textasciicircum 2=(x2,y3)\textasciicircum 2=(x3,y1)\textasciicircum 2= 1,

(x3,y2)\textasciicircum 2=(x3,y3)\textasciicircum 2=1,

((x1,y1),x1)=((x1,y1),x2)=((x1,y1),x3)=((x1,y1),y1)=

((x1,y1),y2)=((x1,y1),y3)=1,

((x1,y2),x1)=((x1,y2),x2)=((x1,y2),x3)=((x1,y2),y1)=

((x1,y2),y2)=((x1,y2),y3)=1,

((x1,y3),x1)=((x1,y3),x2)=((x1,y3),x3)=((x1,y3),y1)=

((x1,y3),y2)=((x1,y3),y3)=1,

((x2,y1),x1)=((x2,y1),x2)=((x2,y1),x3)=((x2,y1),y1)=

((x2,y1),y2)=((x2,y1),y3)=1,

((x2,y2),x1)=((x2,y2),x2)=((x2,y2),x3)=((x2,y2),y1)=

((x2,y2),y2)=((x2,y2),y3)=1,

((x2,y3),x1)=((x2,y3),x2)=((x2,y3),x3)=((x2,y3),y1)=

((x2,y3),y2)=((x2,y3),y3)=1,

((x3,y1),x1)=((x3,y1),x2)=((x3,y1),x3)=((x3,y1),y1)=

((x3,y1),y2)=((x3,y1),y3)=1,

((x3,y2),x1)=((x3,y2),x2)=((x3,y2),x3)=((x3,y2),y1)=

((x3,y2),y2)=((x3,y2),y3)=1,

((x3,y3),x1)=((x3,y3),x2)=((x3,y3),x3)=((x3,y3),y1)=

((x3,y3),y2)=((x3,y3),y3)=1
>;}
\smallskip

Construct the largest 2-quotient group of $\II(3)$ having lower exponent-$2$ class at most 100 as group in the category GrpPC:
\smallskip

{\tt G3,q:=pQuotient(G,2,100);}

Order of $G3$ (The result shows that $|G3|=\II(3)$, and so $G3\cong \II(3)$):\smallskip

{\tt FactoredOrder(G3);}

Construct the graph $\Gamma(3)=\Cay(G3, S)$:
\smallskip

{\tt x1:=x1@q;
x2:=x2@q;
x3:=x3@q;
y1:=y1@q;
y2:=y2@q;
y3:=y3@q;

X:=sub<G3|x1,x2,x3>;
Y:=sub<G3|y1,y2,y3>;

S:={x:x in X|x ne G3!1} join {y:y in Y|y ne G3!1};

Gamma3:=Cay(G3,S);}
\smallskip

Automorphism Group of $\Gamma(3)$ (The result shows that $|\Aut(\Gamma(3))|=|G(3)||\GL_3(2)\wr S_2|$):

{\tt A:=AutomorphismGroup(Gamma3);

$\sharp$A eq ($\sharp$G3)*$\sharp$(GL(3,2))*$\sharp$(GL(3,2))*2;}

\section*{Acknowledgements}

The first author has been supported by the Croatian
Science Foundation under the project 6732. The third author was supported by the National Natural Science Foundation of China (12071023,1211101360).

\end{document}